\input amstex
\NoBlackBoxes

\documentstyle{conm-p}
\hfuzz=10pt

\input xy
\xyoption{all}

%
\def\Z{\Bbb Z} 
\def\Hom{\operatorname{Hom}} 
\def\ker{\operatorname{ker}} 
\def\im{\operatorname{im}} 
\def\coker{\operatorname{coker}} 
\def\t{\thinspace}
\def\di{\displaystyle}
\topmatter
\title ALGEBRAIC POINCAR\'E COBORDISM\endtitle
\author ANDREW RANICKI\endauthor
\leftheadtext{ANDREW RANICKI}%
\rightheadtext{ALGEBRAIC POINCAR\'E COBORDISM}%
\address 
Department of Mathematics and Statistics\newline
\indent University of Edinburgh\newline
\indent Edinburgh EH9 3JZ \newline
\indent Scotland, UK\newline 
~
\endaddress
\email aar\@maths.ed.ac.uk\endemail
\issueinfo{00}
{}
{}
{2000}
\endtopmatter
\document
\noindent {\bf Introduction}
\medskip

The object of this paper is to give a reasonably leisurely account of
the algebraic Poincar\'e cobordism theory of Ranicki [15],\t [16] and the further
development due to Weiss [19], along with some of the applications 
to manifolds and vector bundles. It is a companion paper to Ranicki [17],
which is an introduction to algebraic surgery using forms and formations.

\medskip

Algebraic Poincar\'e cobordism is modelled on the bordism groups
$\Omega_*(X)$ of maps $f:M \to X$ from manifolds to a space $X$.  
The Wall [18] surgery obstruction groups $L_*(A)$ of a ring with involution $A$
were expressed 
in [15] as the cobordism groups of $A$-module chain complexes $C$ with a
quadratic Poincar\'e duality
$$\psi~:~H^{n-*}(C)~\cong~H_*(C)~,$$
and the surgery obstruction $\sigma_*(f,b)\in L_n(\Z[\pi_1(X)])$ of an
$n$-dimensional normal map $(f,b):M\to X$ was expressed as the
cobordism class $(C,\psi)$ of an $n$-dimensional f.g.  free
$\Z[\pi_1(X)]$-module chain complex $C$ such that
$$H_*(C)~=~\ker(f_*:H_*(\widetilde{M})\to H_*(\widetilde{X}))$$
together with an $n$-dimensional quadratic Poincar\'e duality $\psi$.
The passage from the bundle map
$b:\nu_M\to \nu_X$ to $\psi$ used an equivariant chain level version of
the relationship established by Thom between the Wu classes of the
normal bundle $\nu_M$ of a manifold $M$ and the action of the Steenrod
algebra on the Thom class of $\nu_M$.  \medskip

A chain bundle $(C,\gamma)$ over a ring with involution $A$ is an $A$-module
chain complex $C$  together with a Tate $\Z_2$-hypercohomology class
$\gamma \in \widehat{H}^0(\Z_2;C^* \otimes C^*)$.  
The $L$-groups $L^*(C,\gamma)$ of [19] are the cobordism
groups of symmetric Poincar\'e complexes over $A$ with a chain bundle
map to $(C,\gamma)$, which are related to the quadratic $L$-groups by
an exact sequence of abelian groups
$$\dots \to L_n(A) \to L^n(C,\gamma) \to Q_n(C,\gamma) \to L_{n-1}(A) \to
\dots$$
with the $Q$-groups $Q_*(C,\gamma)$ defined purely homologically.  
The surgery obstruction groups $L_*(A)$ of [18] and the symmetric
$L$-groups $L^*(A)$ of Mishchenko [13] are particular examples of the
generalized $L$-groups $L^*(C,\gamma)$.  The main novelty of this paper
is an explicit formula obtained in \S7
for the addition of elements in $Q_*(C,\gamma)$. 
The Wu classes $v_*(\nu) \in H^*(X;\Z_2)$ of a $(k-1)$-spherical
fibration $\nu$ over a space $X$ (e.g. the sphere bundle of a 
$k$-plane bundle) determine a chain bundle $(C(\widetilde{X}),\gamma(\nu))$
over $\Z[\pi_1(X)]$, with $\widetilde{X}$ the universal cover of $X$,
and with a morphism 
$$\pi_{n+k}(T(\nu)) \to Q_n(C(\widetilde{X}),\gamma)~.$$ 
For a $k$-plane bundle $\nu$ the 
morphism factors through the flexible signature map of [19]
$$\Omega_n(X,\nu)~=~\pi_{n+k}(T(\nu)) \to L^n(C(\widetilde{X}),\gamma(\nu))~.$$ 
with $\Omega_n(X,\nu)$ the bordism group of normal maps $(f:M
\to X,b:\nu_M \to \nu)$ from $n$-dimensional manifolds.

\medskip

In subsequent joint work with Frank Connolly a computation of 
$Q_*(C,\gamma)$ will be used to compute the Cappell Unil-groups
in certain special cases.
\medskip

The titles of the sections are
\medskip

{\parskip=2pt
\parindent=75pt
\item{1.} Rings with involution
\item{2.} Chain complexes
\item{3.} Symmetric, quadratic and hyperquadratic structures
\item{4.} Algebraic Wu classes
\item{5.} Algebraic Poincar\'e complexes
\item{6.} Chain bundles
\item{7.} Normal complexes
\item{8.} Normal cobordism
\item{9.} Normal Wu classes
\item{10.} Forms
\item{11.} An example.}

\medskip
\medskip

\noindent \S1. {\bf Rings with involution}
\medskip

In \S1 we show how an involution $a \mapsto \overline{a}$
on a ring $A$ determines a duality involution functor 
$$\hbox{(f.g. projective left $A$-modules)} \to
\hbox{(f.g. projective left $A$-modules)}~.$$
More generally, duality can be defined using an
antistructure on $A$ in the sense of Wall [18], and the $L$-theory results
described in this paper all have versions for rings with antistructure. 
\medskip

Let $A$ be an associative ring with 1, together with an involution 
$$\raise4pt\hbox{$\overline{~}$}~:~A \to A~;~a \mapsto \overline{a}~,$$
that is a function satisfying
$$\overline{a+b}~=~\overline{a}+\overline{b}~,~
\overline{\overline{a}}~=~a~,~\overline{ab}~=~\overline{b}\,.\,\overline{a}~,~
\overline{1}~=~1 \in A~~(a,b \in A)~.$$
In the topological applications $A=\Z[\pi]$ is a group ring, for
some group $\pi$ equipped with a morphism $w:\pi \to \Z_2=\{+1,-1\}$, 
and the involution is defined by
$$\raise4pt\hbox{$\overline{~}$}~:~A \to A~;~\sum\limits_{g \in \pi}n_gg \mapsto 
\sum\limits_{g \in \pi}n_gw(g)g^{-1}~.$$

We take $A$-modules to be left $A$-modules,
unless a right $A$-action is expressly specified. 
Given an $A$-module $M$ there is defined a right $A$-module
$M^t$ with the same additive group and 
$$M^t \times A \to M^t~;~(x,a) \mapsto \overline{a}x~.$$
The {\it dual} of an $A$-module $M$ is the $A$-module with additive group 
$$M^*~=~\Hom_A(M,A)$$
and $A$ acting by 
$$A \times M^* \to M^*~;~(a,f) \mapsto (x \mapsto f(x).\overline{a})~.$$
The {\it dual} of an $A$-module morphism $f \in \Hom_A(M,N)$ 
is the $A$-module morphism
$$f^*~:~N^* \to M^*~;~g \mapsto (x \mapsto g(f(x)))~.$$ 
For a f.g. (finitely generated) projective $A$-module $M$ 
the natural $A$-module isomorphism 
$$M \to M^{**}~;~x \mapsto (f \mapsto \overline{f(x)})$$ 
will be used to identify 
$$M^{**}~=~M~.$$

\noindent \S2. {\bf Chain complexes} 
\medskip

In order to adequately deal with the quadratic nature of the function
$C\to  C^t\otimes _AC$ sending an
$A$-module chain complex $C$ to the $\Z$-module chain complex 
$C^t\otimes_AC$ it is necessary to use the equivalence of Dold [7]
and Kan [8] (cf. May [10,\t\S 22])
between positive $\Z$-module chain complexes and simplicial $\Z$-modules. 
This is now recalled, along with some other properties of chain complexes that
we shall require. 
\medskip

An $A$-module chain complex 
$$C~ :~ \dots\to C_{r+1}~\raise4pt\hbox{$d \atop \to$}~C_r~
\raise4pt\hbox{$d \atop \to$}~C_{r-1} \to \dots~~ (r\in \Z)$$
is {\it $n$-dimensional} if each $C_r$ ($0 \le r \le n$)
is a f.g. projective $A$-module and $C_r=0$ for $r<0$ and $r>n$. 
By an abuse of terminology chain
complexes of the chain homotopy type of an $n$-dimensional chain complex
will also be called $n$-dimensional. 
\medskip

The {\it suspension} of an
$A$-module chain complex $C$ is the $A$-module chain complex defined by
$$d_{SC}~=~d_C~ :~ SC_r~=~C_{r-1} \to SC_{r-1}~=~C_{r-2}~.$$
If $C$ is $n$-dimensional then $SC$ is $(n+1)$-dimensional. 
\medskip

Given an $A$-module chain complex $C$ let 
$$C^r~=~(C_r)^*~~(r\in \Z) ~.$$ 
The {\it dual $A$-module} chain complex $C^{-*}$ is defined by 
$$d_{C^{-*}}~=~(d_C)^*~:~(C^{-*})_r~=~C^{-r} \to 
(C^{-*})_{r-1}~=~C^{-r+1}~.$$ 
The {\it $n$-dual $A$-module chain complex} $C^{n-*}$ is defined by
$$d_{C^{n-*}}~=~(-1)^r(d_C)^*~:~
(C^{n-*})_r~=~C^{n-r} \to (C^{n-*})_{r-1}~=~C^{n-r+1}~.$$
The $n$-fold suspension of the dual $S^nC^{-*}$ is related to the $n$-dual
$C^{n-*}$ by the isomorphism 
$$S^nC^{-*} \to C^{n-*}~;~x \mapsto (-1)^{r(r-1)/2}x~ (x\in C^{n-r})~.$$
In particular, $C^{0-*}$ is isomorphic to (but not identical to) $C^{-*}$.
\medskip

A {\it chain map up to sign} between $A$-module chain complexes
$$ f~:~C \to D$$
is a collection of $A$-module morphisms 
$$\{f\in \Hom_A(C_r,D_r)\,\vert\, r\in \Z\}$$ 
such that
$$d_Df~=~\pm fd_C~ :~ C_r \to D_{r-1}~~ (r\in \Z)~ .$$
If the sign is always + this is just a {\it chain map} $f:C\to D$, as
usual. 
\medskip

Given $A$-module chain complexes $C,D$ let $C^t\otimes_AD$,
$\Hom_A(C,D)$ be the $\Z$-module chain complexes defined by
$$\eqalign{
&(C^t\otimes _AD)_n~=~\sum\limits_{p+q=n}C_p \otimes _AD_q~,\cr
&\hskip100pt d_{C^t\otimes _AD}(x\otimes y)~=~x\otimes d_D(y) +
(-1)^qd_C(x)\otimes y~,\cr
&\Hom_A(C,D)_n~=~\sum\limits_{q-p=n}\Hom_A(C_{p},D_q)~,\cr
&\hskip100pt d_{\Hom_A(C,D)}(f)(x)~=~d_D(f(x)) + (-1)^qf(d_C(x))~.}$$
A cycle $f\in \Hom_A(C,D)_n$ is a chain map up to sign $f:S^nC\to D$, and 
$$H_n(\Hom_A(C,D))~=~H_0(\Hom_A(S^nC,D))$$
is the $\Z$-module of chain homotopy classes of chain maps $S^nC\to D$. 
\medskip

For finite-dimensional $C$ the {\it slant isomorphism} of $\Z$-module chain
complexes 
$$C^t\otimes _AD \to \Hom_A(C^{-*},D)~ ;~ 
x\otimes y \mapsto (f \mapsto \overline{f(x)}\,.\,y)$$
will be used to identify 
$$C^t\otimes _AD~=~\Hom_A(C^{-*},D)~ .$$
A cycle 
$$f \in (C^t\otimes _AD)_n~=~\Hom_A(C^{-*},D)_n$$ 
is a chain map $f : C^{n-*}\to D$. Thus 
$$H_n(C^t\otimes_AD)~=~H_n(\Hom_A(C^{-*},D))~=~H_0(\Hom_A(C^{n-*},D))$$
is the $\Z$-module of chain homotopy classes of chain maps $C^{n-*}\to D$. 
\medskip

The {\it algebraic mapping cone} $C(f)$ of an $A$-module chain map
$f:C\to D$ is the $A$-module chain complex defined as usual by
$$d_{C(f)}~=~\pmatrix d_D& (-1)^{r-1}f \cr 0& d_C\endpmatrix~ :~
C(f)_r~=~D_r\oplus C_{r-1} \to C(f)_{r-1}~=~D_{r-1}\oplus C_{r-2}~.$$
The relative homology $A$-modules 
$$H_n(f)~=~H_n(C(f))$$ 
are such that there is defined an exact sequence 
$$\dots \to H_n(C)~\raise4pt\hbox{$f_* \atop \to$}~H_n(D) \to 
H_n(f) \to H_{n-1}(C) \to \dots~ .$$

Let $C(\Delta^n)$ denote the cellular chain complex of the standard
$n$-simplex $\Delta^n$ with the standard cell structure consisting of
$\displaystyle{n+1 \choose r+1}$ $r$-cells $(0 \le r \le n)$.
\medskip

Given a $\Z$-module chain
complex $C$ let $K(C)$ denote the simplicial $\Z$-module defined by the
Dold-Kan construction, with one $n$-simplex for each chain map
$C(\Delta^n)\to C$ and the evident face and degeneracy maps
$d_i,s_i$, such that 
$$\pi_n(K(C))~=~H_n(C)~~(n\ge 0)~.$$
Given a chain $y\in C_n$ and cycles
$x_i\in \ker(d:C_{n-1}\to C_{n-2})$ $(0 \le i \le n)$ such that
$$dy~=~\sum\limits^n_{i=0}(-1)^ix_i \in C_{n-1}$$ 
let $(y;x_0,\dots,x_n)$ denote the $n$-simplex of $K(C)$ defined by the
chain map 
$$f~ :~ C(\Delta^n) \to C$$ 
with 
$$\eqalign{
&f~:~C(\Delta^n)_n~=~\Z \to C_n~ ;~ 1 \mapsto y~,\cr
&d_if~ :~C(\Delta^{n-1})_{n-1}~=~\Z \to C_{n-1}~;~1 \mapsto x_i~~
(0 \le i \le n)~.}$$
The chain $x\in C_n$ is
identified with the $n$-simplex $(x;dx,0,\dots,0)\in K(C)^{(n)}$. 
\medskip

Given $\Z$-module chain complexes $C,D$ and a simplicial map 
$$f~ :~ K(C) \to K(D)$$ 
(which need not preserve the $\Z$-module structure) there is defined
a function 
$$f~ : ~C_n \to D_n~;~x \mapsto f(x)~=~f(x;dx,0,\dots,0)$$
such that
$$df(x)~=~f(dx) \in D_{n-1}~~(x \in C_n)~.$$
In general, $f:C_n \to D_n$ is only additive on the level of homology, with
$$f(x+x') -f(x) - f(x')~=~d[x,x']_f \in D_n~~(x,x' \in\ker(d:C_n\to C_{n-1})$$ 
where the function
$$[~,~]_f~:~\ker(d:C_n\to C_{n-1}) \times \ker(d:C_n\to C_{n-1})\to D_{n+1}$$
is defined by
$$[x,x']_f~=~f(0;x+x',x,-x',0,\dots,0)-f(0;x,x,0,\dots,0)-f(0;x',0,-x',0,\dots,0)~.$$
The induced functions
$$f_*~:~H_n(C) \to H_n(D)~;~x \mapsto f(x)$$
are thus morphisms of abelian groups, which fit into an exact sequence 
$$\dots \to H_{n+1}(f) \to H_n(C)~
\raise4pt\hbox{$f_* \atop\to$}~H_n(D) \to H_n(f) \to  H_{n-1}(C) \to \dots~,$$
with the relative group $H_n(f)$ ($= \pi_n(f)$) the set of equivalence classes
of pairs $(x,y) \in C_n\times D_{n+1}$ such that 
$$dx~=~0 \in C_{n-1}~,~f(x)~=~dy \in D_n~,$$ 
subject to the equivalence relation
$$\eqalign{
(x,y) \sim (x',y')~~&\hbox{if there exist}~(u,v) \in 
C_{n+1}\times D_{n+2}~\hbox{such that}\cr
&x - x'~=~du \in C_n~,~y-y'~=~f(u;x,x',0,\dots,0) + dv \in D_{n+1}~,}$$ 
and addition by
$$(x,y) + (x',y')~=~(x+x',y+y'+[x,x']_f) \in H_n(f)~.$$ 
If $f:K(C)\to K(D)$ does preserve the $\Z$-module structure (so that
$[~,~]_f=0$) then $f$ is essentially just a chain map $f:C\to D$, and
the relative homology groups $H_*(f)$ are just the homology groups
$H_*(C(f))$ of the algebraic mapping cone $C(f)$, as usual.  
\medskip

\noindent \S3. {\bf Symmetric, quadratic and hyperquadratic structures} 
\medskip

An $n$-dimensional 
$\cases \text{symmetric} \cr \text{quadratic} \cr \text{hyperquadratic}\endcases$
structure on an $A$-module chain complex $C$ is a cycle 
representing an element of the
$\cases \text{$\Z_2$-hypercohomology}\cr
\text{$\Z_2$-hyperhomology}\cr
\text{Tate $\Z_2$-hypercohomology}\endcases$ 
group
$\cases H^n(\Z_2;C^t\otimes_AC)\cr
H_n(\Z_2;C^t\otimes_AC)\cr
\widehat{H}^n(\Z_2;C^t\otimes_AC)\endcases$
in the sense of Cartan and Eilenberg [6].
\medskip

Given an $A$-module chain complex $C$
let the generator $T\in \Z_2$ act on $C^t\otimes _AC$ by the
{\it transposition involution}
$$T~ :~ C^t\otimes _AC \to C^t\otimes _AC~ ;~ 
x\otimes y~\mapsto ~(-1)^{pq}y\otimes x ~~ (x\in C_p,y\in C_q)~ .$$
For finite-dimensional $C$ use the slant isomorphism to identify
$$C^t\otimes _AC~=~\Hom_A(C^{-*},C)~ .$$
Under this identification the transposition involution corresponds to
the duality involution on $\Hom_A(C^{-*},C)$ 
$$T~:~\Hom_A(C^{-*},C) \to \Hom_A(C^{-*},C)~;~\phi\mapsto (-1)^{pq}\phi^*~
(\phi\in \Hom_A(C^p,C_q))~.$$
A cycle $\phi\in \Hom_A(C^{-*},C)_n=(C^t\otimes _AC)_n$ 
is a chain map $\phi:C^{n-*}\to C$, and $H_n(\Hom_A(C^{-*},C))$ is the
$\Z$-module of chain homotopy classes of $A$-module chain maps
$C^{n-*}\to C$. 
Let $W$ be the standard free $\Z[\Z_2]$-module resolution of $\Z$ 
$$W~:~\dots \to \Z[\Z_2]~
\raise4pt\hbox{$1-T \atop \to$}~\Z[\Z_2]~
\raise4pt\hbox{$1+T \atop \to$}~\Z[\Z_2]~
\raise4pt\hbox{$1-T \atop \to$}~\Z[\Z_2] \to 0$$
and let $\widehat{W}$ be the complete resolution 
$$\widehat{W}~:~\dots \to \Z[\Z_2]~
\raise4pt\hbox{$1-T \atop \to$}~\Z[\Z_2]~
\raise4pt\hbox{$1+T \atop \to$}~\Z[\Z_2]~
\raise4pt\hbox{$1-T \atop \to$}~\Z[\Z_2] \to \dots~.$$
The $\cases \hbox{\it $\Z_2$-hypercohomology}\cr
\hbox{\it $\Z_2$-hyperhomology}\cr
\hbox{\it Tate $\Z_2$-hypercohomology}\endcases$
{\it groups} of a $\Z[\Z_2]$-module chain complex $C$ are defined by
$$\cases H^n(\Z_2;C)~=~H_n(\Hom_{\Z[\Z_2]}(W,C))\cr
H_n(\Z_2;C)~=~H_n(W\otimes_{\Z[\Z_2]}C)\cr
\widehat{H}^n(\Z_2;C)~=~H_n(\Hom_{\Z[\Z_2]}(\widehat{W},C))~.\endcases$$
The evident short exact sequence of $\Z[\Z_2]$-module chain complexes
$$0 \to S W^{-*} \to \widehat{W} \to W \to 0$$
induces a long exact sequence of abelian groups 
$$\dots\to H_n(\Z_2;C)~\raise4pt\hbox{$1+T \atop \to$}~H^n(\Z_2;C)
\raise4pt\hbox{$J \atop \to$}~\widehat{H}^n(\Z_2;C)
\raise4pt\hbox{$H \atop \to$}~H_{n-1}(\Z_2;C) \to \dots~.$$
An element 
$\cases \phi\in H^n(\Z_2;C)\cr
\psi\in H_n(\Z_2;C)\cr
\theta\in \widehat{H}^n(\Z_2;C)\endcases$ is
represented by an $n$-cycle of 
$\cases \Hom_{\Z[\Z_2]}(W,C)\cr
W\otimes_{\Z[\Z_2]}C\cr
\Hom_{\Z[\Z_2]}(\widehat{W},C)\endcases$
which is just a collection of chains of $C$
$\cases 
\{\phi_s\in C_{n+s}\,\vert\, s\ge 0 \}\cr
\{\psi_s\in C_{n-s}\,\vert\, s \ge 0 \}\cr
\{\theta_s\in C_{n+s}\,\vert\, s \in \Z\}\endcases$ 
such that 
$$\cases 
d_C(\phi_s) +
(-1)^{n+s-1}(\phi_{s-1}+(-1)^sT\phi_{s-1})~=~0 \in C_{n+s-1}~(s\ge 0,\phi_{-1}=0)\cr
d_C(\psi_s) +(-1)^{n-s-1}(\psi_{s+1}+(-1)^{s+1}T\psi_{s+1})~=~0 \in C_{n-s-1}~
(s \ge 0)\cr
d_C(\theta_s) +(-1)^{n+s-1}(\theta_{s-1}+(-1)^sT\theta_{s-1})~=~0 
\in C_{n+s-1}~(s\in\Z)\endcases$$
with 
$$\eqalign{
&1+T~:~H_n(\Z_2;C) \to H^n(\Z_2;C)~ ;\cr
&\hskip50pt \psi~=~\{\psi_s\,\vert\, s \ge 0\} \mapsto \{
(1+T)\psi_s=\cases (1+T)\psi_0&\text{if $s=0$}\cr
0&\text{if $s\ge 1$}\endcases\}~,\cr
&J~:~H^n(\Z_2;C) \to \widehat{H}^n(\Z_2;C)~ ;\cr
&\hskip50pt \phi~=~\{\phi_s\,\vert\, s\ge 0\} \mapsto 
\{J\phi_s=\cases \phi_s&\text{if $s \ge 0$} \cr
0&\text{if $s \le -1$}\endcases\}~,\cr
&H~:~\widehat{H}^n(\Z_2;C) \to H_{n-1}(\Z_2;C)~ ;\cr
&\hskip50pt \theta~=~\{\theta_s \,\vert\, s \in \Z\} \mapsto 
 H\theta~=~\{H\theta_s=\theta_{-s-1}\,\vert\, s \ge 0\}~.}$$
Given an $A$-module chain complex $C$ use the action of 
$T \in \Z_2$ on $C^t\otimes _AC$ by the transposition involution to 
define the $\Z$-module chain complex 
$$\eqalign{
&W^{\%}C~=~\Hom_{\Z[\Z_2]}(W,C^t\otimes _AC)\cr
&W_{\%}C~=~W\otimes_{\Z[\Z_2]}(C^t\otimes _AC)\cr
&\widehat{W}^{\%}C~=~\Hom_{\Z[\Z_2]}(\widehat{W},C^t\otimes _AC)~.}$$
We shall be mainly concerned with finite-dimensional $C$, using the slant 
isomorphism to identify 
$$C^t\otimes _AC~=~\Hom_A(C^{-*},C)$$ 
and 
$$\eqalign{
&W^{\%}C~=~\Hom_{\Z[\Z_2]}(W,\Hom_A(C^{-*},C))\cr
&W_{\%}C~=~W\otimes_{\Z[\Z_2]}\Hom_A(C^{-*},C)\cr
&\widehat{W}^{\%}C~=~\Hom_{\Z[\Z_2]}(\widehat{W},\Hom_A(C^{-*},C))~.}$$
An {\it $n$-dimensional $\cases 
\text{\it symmetric} \cr 
\text{\it quadratic} \cr
\text{\it hyperquadratic}\endcases$ structure}
on a finite-dimensional $A$-module chain complex $C$ is a cycle 
$\cases \phi\in (W^{\%}C)_n\cr 
\psi\in (W_{\%}C)_n\cr
\theta\in (\widehat{W}^{\%}C)_n~,\endcases$
which is just a collection of $A$-module morphisms 
$$\cases
\{\phi_s\in \Hom_A(C^{n-r+s},C_r)\,\vert\, r\in \Z,s \ge 0 \}\cr
\{\psi_s\in \Hom_A(C^{n-r-s},C_r)\,\vert\, r\in \Z,s \ge 0 \}\cr
\{\theta_s\in \Hom_A(C^{n-r+s},C_r)\,\vert\, r\in \Z,s \in \Z \}
\endcases$$
such that 
$$\cases
d\phi_s+(-1)^r\phi_sd^*+(-1)^{n+s-1}(\phi_{s-1}+(-1)^sT\phi_{s-1})~=~0\cr
\hskip150pt :~C^{n-r+s-1} \to C_r~~ (s\ge 0,\phi_{-1}=0)\cr
d\psi_s+(-1)^r\psi_sd^*+(-1)^{n-s-1}(\psi_{s+1}+(-1)^{s+1}T\psi_{s+1})~=~0\cr
\hskip150pt :~ C^{n-r-s-1} \to C_r~~ (s\ge 0)\cr
d\theta_s+(-1)^r\theta_sd^*+(-1)^{n+s-1}(\theta_{s-1}+(-1)^sT\theta_{s-1})~=~0\cr
\hskip150pt :~C^{n-r+s-1} \to C_r~~ (s\in \Z)~.
\endcases$$

An {\it equivalence} $\cases \xi:\phi \to \phi'\cr
\chi:\psi\to \psi'\cr
\nu:\theta \to \theta'\endcases$ of $n$-dimensional
$\cases \text{symmetric} \cr 
\text{quadratic} \cr 
\text{hyperquadratic}\endcases$
structures on $C$ is a chain 
$\cases \xi\in (W^{\%}C)_{n+1}\cr
\chi\in (W_{\%}C)_{n+1}\cr
\nu\in (\widehat{W}^{\%}C)_{n+1}\endcases$ such that 
$$\cases \phi' - \phi~=~d(\xi) \in (W^{\%}C)_n\cr
\psi' - \psi~=~d(\chi) \in (W_{\%}C)_n\cr
\theta' - \theta~=~d(\nu) \in (\widehat{W}^{\%}C)_n~.\endcases$$
The {\it $n$-dimensional 
$\cases \text{\it symmetric}\cr \text{\it quadratic}\cr 
\text{hyperquadratic}\endcases$ structure group}
$\cases Q^n(C) \cr Q_n(C) \cr \widehat{Q}^n(C)\endcases$
of a chain complex $C$ is the abelian group
of equivalence classes of $n$-dimensional
$\cases \text{symmetric}\cr \text{quadratic}\cr \text{hyperquadratic}\endcases$
structures on $C$, that is 
$$\cases 
Q^n(C)~=~H^n(\Z_2;C^t\otimes_AC)~=~H_n(W^{\%}C)\cr
Q_n(C)~=~H_n(\Z_2;C^t\otimes_AC)~=~H_n(W_{\%}C)\cr
\widehat{Q}^n(C)~=~\widehat{H}^n(\Z_2;C^t\otimes_AC)~=~
H_n(\widehat{W}^{\%}C)~.\endcases$$

The $Q$-groups are related by a long exact sequence 
$$\dots \to Q_n(C)~
\raise4pt\hbox{$1+T \atop \to$}~Q^n(C)~
\raise4pt\hbox{$J \atop \to$}~\widehat{Q}^n(C)~
\raise4pt\hbox{$H \atop \to$}~Q_{n-1}(C) \to \dots$$
involving the morphisms induced in
homology by the $\Z$-module chain maps 
$$1+T~:~ W_{\% }C \to W^{\% }C~ ,~ J~ :~
W^{\% }C \to \widehat{W}^{\%}C~,~H~:~\widehat{W}^{\%}C \to 
S(W_{\%}C)$$
defined by 
$$\eqalign{
&1+T~:~(W_{\%}C)_n \to (W^{\%}C)_n~ ;\cr
&\hskip25pt \{\psi_s\in (C^t\otimes_AC)_{n-s}\,\vert\, s \ge 0\} \mapsto
\{((1+T)\psi)_s=\cases (1+T)\psi_0&\text{if $s=0$}\cr 
0&\text{if $s \ge 1$}\endcases \}~,\cr
&J~ :~ (W^{\%}C)_n \to (\widehat{W}^{\%}C)_n~ ;\cr
&\hskip25pt \{\phi_s\in (C^t\otimes_AC)_{n+s}\,\vert\, s \ge 0\} \mapsto 
\{(J\phi)_s=
\cases \phi_s&\text{if $s \ge 0$} \cr 0&\text{if $s \le -1$}
\endcases\}~,\cr
&H~:~(\widehat{W}^{\%}C)_n \to (W_{\% }C)_{n-1}~ ;\cr
&\hskip25pt \{\theta_s\in(C^t\otimes _AC)_{n+s}\,\vert\,s\in\Z\}~
\mapsto~\{(H\theta)_s=\theta_{-s-1}\,\vert\, s \ge 0\}~.}$$
An $n$-dimensional symmetric structure $\phi\in (W^{\% }C)_n$ is equivalent
to the symmetrization $(1+T)\psi$ of an $n$-dimensional quadratic
structure $\psi\in (W_{\%}C)_n$ if and only if the $n$-dimensional
hyperquadratic structure $J(\phi)\in (\widehat{W}^{\% }C)_n$ is
equivalent to 0. 
An $A$-module chain map $f:C\to  D$ induces a $\Z[\Z_2]$-module chain map
$$f^t\otimes _Af~:~C^t\otimes _AC \to D^t\otimes _AD~ ;~
 x\otimes y \mapsto f(x)\otimes f(y)$$
and hence $\Z$-module chain maps
$$\eqalign{&f^{\%}~:~W^{\%}C \to W^{\%}D~,\cr
&f_{\%}~:~W_{\%}C \to W_{\%}D~,\cr
&\widehat{f}^{\%}~:~\widehat{W}^{\%}C \to \widehat{W}^{\%}D~.}$$
An $A$-module chain homotopy
$$g~ :~ f \simeq f'~ :~ C \to D$$
determines $\Z$-module chain homotopies
$$\eqalign{&
(g;f,f')^{\%}~:~f^{\%}\simeq f^{\prime \%}~:~W^{\%}C \to W^{\%}D\cr
&(g;f,f')_{\%}~:~f_{\%}\simeq f'_{\%}~~:~W_{\%}C \to W_{\%}D\cr
&\widehat{(g;f,f')}^{\%}~:~\widehat{f}^{\%}\simeq \widehat{f}^{\prime \%}~:~
\widehat{W}^{\%}C \to \widehat{W}^{\%}D}$$
with
$$\eqalign{&
(g;f,f')^{\% }~:~(W^{\% }C)_n~=~\sum\limits^{\infty}_{s=0}
(C^t\otimes _AC)_{n+s}\cr
&\hskip100pt \to (W^{\% }D)_{n+1}~=~\sum\limits^{\infty}_{s=0}\sum\limits_q
D^t_{n-q+s+1}\otimes _AD_q~;\cr
&\hskip50pt \sum\limits^{\infty}_{s=0}\phi_s \mapsto 
\sum\limits^{\infty}_{s=0}((f^t\otimes _Ag + g^t\otimes _Af')(\phi_s) + 
(-1)^{q+s-1}(g^t\otimes _Ag)(T\phi_{s-1}))}$$
and similarly for
$(g;f,f')_{\% }$, $\widehat{(g;f,f')}^{\% }$. 
Thus the induced morphisms in the $Q$-groups
$$\eqalign{&f^{\%}~:~Q^n(C) \to Q^n(D)\cr
&f_{\%}~:~Q_n(C) \to Q_n(D)\cr
&\widehat{f}^{\%}~:~\widehat{Q}^n(C) \to \widehat{Q}^n(D)}$$
depend only on the chain homotopy class of $f$, 
and are isomorphisms if $f$ is a chain equivalence. 
For finite-dimensional $C,D$ the slant isomorphisms are
used to identify 
$f^t\otimes _Af:C^t\otimes _AC\to D^t\otimes _AD$ with
$$\Hom_A(f^*,f)~ :~ \Hom_A(C^{-*},C) \to \Hom_A(D^*,D)~ ;~
\theta \mapsto f\theta f^*~,$$
and similarly for
$f^{\%},f_{\%},\widehat{f}^{\%}$ and $(g;f,f')^{\%},(g;f,f')_{\%},
\widehat{(g;f,f')}^{\%}$. 
\medskip

Although all the $Q$-groups are chain homotopy invariant, only
the hyperquadratic $Q$-groups $\widehat{Q}^*(C)$ are additive. The sum
of $A$-module chain maps $f,g:C\to D$ is an
$A$-module chain map $f+g : C\to D$
such that
$$\eqalign{
&(f+g)^{\%} - f^{\%} - g^{\%}~:~ Q^n(C) \to H_n(C^t\otimes_AC)~
\raise4pt\hbox{$f^t\otimes _Ag\atop \to$}~H_n(D^t\otimes_AD)\to Q^n(D)~,\cr
&(f+g)_{\%} - f_{\%} - g_{\%}~:~ Q_n(C) \to H_n(C^t\otimes_AC)~
\raise4pt\hbox{$f^t\otimes _Ag\atop \to$}~H_n(D^t\otimes_AD)\to Q_n(D)~,\cr
&\widehat{(f+g)}^{\%} -\widehat{f}^{\%} - \widehat{g}^{\%}~=~0~:~
\widehat{Q}^n(C) \to \widehat{Q}^n(D)}$$
with 
$$\eqalign{
&Q^n(C) \to H_n(C^t\otimes_AC)~ ;~
\phi~=~\{\phi_s\,\vert\, s \ge 0\} \mapsto \phi_0~,\cr
&Q_n(C) \to H_n(C^t\otimes_AC)~;~
\psi~=~\{\psi_s\,\vert\, s \ge 0\} \mapsto (1+T)\psi_0~,\cr
&H_n(D^t\otimes _AD) \to Q^n(D)~;~
\theta \mapsto \{\phi_s=\cases (1+T)\theta&\text{if $s=0$}\cr 
0&\text{if $s \ge 1$}\endcases \}~,\cr
&H_n(D^t\otimes _AD) \to Q_n(D)~;~
\theta \mapsto \{\psi_s=\cases \theta&\text{if $s=0$} \cr 0&\text{if $s \ge 1$}\endcases\}~.}$$
Given a finite-dimensional $A$-module chain complex $C$ and $n \ge 0$
define the {\it $n$-fold suspension chain isomorphism} 
$$\eqalign{
S^n~:~&S^n(\widehat{W}^{\%}C) \to \widehat{W}^{\%}(S^nC)~;\cr
&\theta~=~\{\theta_s\in \Hom_A(C^r,C_{m-r+s})\,\vert\, s\in \Z\}\cr
&\mapsto~S^n\theta~=~\{(S^n\theta)_s=\theta_{s-n}\in 
\Hom_A(C^r,C_{m-n+r+s})\,\vert\, s\in \Z\}~.}$$
For any (finite-dimensional) $A$-module chain
complexes $C,D$ there is defined a simplicial map 
$$I~ :~ K(\Hom_A(C,D)) \to K(\Hom_\Z(\widehat{W}^{\%}C,
\widehat{W}^{\%}D))$$
sending a cycle $f\in \Hom_A(C,D)_n$ 
(= a chain map up to sign $f:S^nC\to D$) to the $\Z$-module chain map up to
sign 
$$I(f)~=~\widehat{f}^{\%}S^n~:~S^n(\widehat{W}^{\%}C)~
\raise4pt\hbox{$S^n\atop \to$}~\widehat{W}^{\%}S^nC~
\raise4pt\hbox{$\widehat{f}^{\%} \atop \to$}~W^{\% }D~ .$$
An $n$-simplex $(g;f,f',0,\dots ,0)\in K(\Hom_A(C,D))^{(n)}$ 
(= an $A$-module chain homotopy up to sign $g:f\simeq f':S^nC\to D$) 
is sent to the $\Z$-module chain homotopy up to sign 
$$I(g;f,f')~=~(g;f,f')^{\% }S^n~:~I(f)~\simeq~I(f')~:~
S^n(\widehat{W}^{\% }C) \to \widehat{W}^{\% }D~.$$ 
The failure of $I$ to be linear on chain maps up to sign $f:S^nC\to D$ 
is given by the chain homotopy up to sign
$$ [f,f']~ :~ \widehat{(f+f')}^{\%}~ \simeq~
\widehat{f}^{\% } + \widehat{f}^{\prime \%} ~:~
\widehat{W}^{\%}(S^nC) \to \widehat{W}^{\%}D$$ 
defined by 
$$\eqalign{
&[f,f']~:~(S^n\widehat{W}^{\%}C)_m \to (\widehat{W}^{\%}D)_{m+n+1}~;\cr
&\theta~=~\{\theta_s\in \Hom_A(C^r,C_{m-r+s})\,\vert\, s\in\Z\}\cr
&\mapsto [f,f']\theta~=~\{T^{n+1}f\theta_{s-n+1}f^{\prime *}\in 
\Hom_A(D^r,D_{m-n+r+s+1})\,\vert\, s\in \Z\}~.}$$

\noindent \S4. {\bf Algebraic Wu classes}
\medskip

The algebraic Wu classes are the fundamental invariants of a duality
structure on a chain complex $C$, which are obtained by an algebraic
analogue of the Steenrod squares in the cohomology groups of a
topological space.  In the topological applications the algebraic Wu
classes are closely related to the topological Wu classes, 
as explained in Ranicki [15].
\medskip

Let $S^rA$ $(r\in \Z)$ denote the $A$-module chain complex 
$$S^rA~ :~ \dots \to 0 \to  A \to  0 \to \dots$$
concentrated in degree $r$. For any $A$-module chain complex $C$ there are
defined natural isomorphisms 
$$H_0(\Hom_A(C,S^rA)) \to H^r(C)~;~(f:C_r\to  A) \mapsto f^*(1)~.$$ 
An element $f\in H^r(C)$ is just a chain homotopy class of chain maps
$f:C\to S^rA$.  The Wu classes of a quadratic structure on $C$ are the
invariants of the equivalence class defined by sending an element $f\in
H^r(C)$ to the induced equivalence class of quadratic structures on
$S^rA$.  The quadratic structure groups of the elementary complexes
$S^rA$ are identified with subquotients of the ground ring $A$.
\medskip

An {\it $A$-group} $M$ is an abelian group together with an $A$-action 
$$A \times M \to M~ ;~ (a,x) \to ax$$ 
such that 
$$a(x+y)~=~ax + ay~~ ,~~ a(bx)~=~(ab)x~~ ,~~ 1x~=~x~~(x,y\in M, a,b\in M)~.$$ 
An $A$-module is an $A$-group $M$ such that also
$$(a+b)x~=~ax + bx \in M~.$$ 
An {\it $A$-morphism} of $A$-groups is a
morphism of abelian groups 
$$f~ :~ M \to N$$
such that 
$$f(ax)~=~af(x) \in N~~ (x\in M,a\in A)~ .$$ 
The set of $A$-group morphisms $f:M\to N$ defines an
abelian group $\Hom_A(M,N)$, with addition by 
$$(f+g)(x)~=~f(x) + g(x) \in N~.$$ 
For $A$-modules $M,N$ the $A$-morphisms $f:M\to N$ coincide with
$A$-module morphisms. 
\medskip

For $\epsilon= \pm 1 $ let the generator $T\in \Z_2$ act on $A$ by the 
$\epsilon$-involution 
$$T_{\epsilon}~:~A \to A~ ;~ a \mapsto \epsilon \overline{a}~.$$
Define the 
$\cases \text{\it $\Z_2$-cohomology}\cr
\text{\it $\Z_2$-homology}\cr
\text{\it Tate $\Z_2$-cohomology}\endcases$ {\it $A$-groups}
$\cases H^*(\Z_2;A,\epsilon) \cr
H_*(\Z_2;A,\epsilon) \cr
\widehat{H}^*(\Z_2;A,\epsilon)\endcases$ by 
$$\cases H^r(\Z_2;A,\epsilon)~=~
\cases
\ker(1-T_{\epsilon}:A \to A)&\text{if $r=0$}\cr
\widehat{H}^r(\Z_2;A,\epsilon)&\text{if $r \ge 1$}\cr
0&\text{if $r<0$}\endcases\cr
H_r(\Z_2;A,\epsilon)~=~\cases
\coker(1-T_{\epsilon}:A\to A)&\text{if $r=0$}\cr
\widehat{H}^{r+1}(\Z_2;A,\epsilon)&\text{if $r \ge 1$}\cr
0&\text{if $r<0$}\endcases\cr
\widehat{H}^r(\Z_2;A,\epsilon)~=~
\ker(1-(-1)^rT_{\epsilon}:A\to A)/
\im(1+(-1)^rT_{\epsilon}:A\to A)~~ (r\in \Z)~.\endcases$$
The $A$-action 
$$A \times \widehat{H}^r(\Z_2;A,\epsilon) \to 
\widehat{H}^r(\Z_2;A,\epsilon)~ ;~ (a,x) \mapsto a x \overline{a}$$ 
defines an $A$-module structure on $\widehat{H}^r(\Z_2;A,\epsilon)$. 
The $A$-actions 
$$\eqalign{&A \times 
H^0(\Z_2;A,\epsilon) \to H^0(\Z_2;A,\epsilon)~ ;~ 
(a,x) \mapsto a x \overline{a}\cr
&A \times H_0(\Z_2;A,\epsilon) \to H_0(\Z_2;A,\epsilon) ~ ;~ 
(a,x) \mapsto a x \overline{a}}$$
are not linear in $A$, and so do not define $A$-module structures. 

For $\epsilon = +1$ the groups
$\cases H^*(\Z_2;A,\epsilon)\cr H_*(\Z_2;A,\epsilon) \cr
\widehat{H}^*(\Z TheA,\epsilon)\endcases$
are denoted by
$\cases H^*(\Z_2;A)\cr H_*(\Z_2;A) \cr
\widehat{H}^*(\Z_2;A).\endcases$
\medskip

The natural $\Z$-module isomorphisms 
$$\eqalign{
&Q^n(S^rA) \to H^{2r-n}(\Z_2;A,(-1)^r)~;~\phi \mapsto \phi_{2r-n}(1)(1)\cr
&Q_n(S^rA) \to H_{n-2r}(\Z_2;A,(-1)^r)~;~\psi \mapsto \psi_{n-2r}(1)(1)\cr
&\widehat{Q}^n(S^rA) \to \widehat{H}^{r-n}(\Z_2;A,(-1)^r)~;~
\theta \mapsto \theta_{2r-n}(1)(1)}$$
will be used as identifications. 
\medskip

The {\it Wu classes} of a
$\cases \text{symmetric}\cr \text{quadratic}\cr \text{hyperquadratic}\endcases$
structure $\cases \phi\in (W^{\%}C)_n\cr
\psi\in (W_{\%}C)_n\cr
\theta\in (\widehat{W}^{\%}C)_n\endcases$ are the invariants of the equivalence 
class of structures defined by the $A$-morphisms 
$$\cases 
v_r(\phi)~:~H^{n-r}(C)~=~H_0(\Hom_A(C,S^{n-r}A))\cr
\hskip50pt \to Q^n(S^{n-r}A)~=~H^{n-2r}(\Z_2;A,(-1)^{n-r})~;~
f \mapsto (f\otimes f)(\phi_{n-2r})\cr
v^r(\psi)~:~H^{n-r}(C)~=~H_0(\Hom_A(C,S^{n-r}A))\cr
\hskip50pt\to Q_n(S^{n-r}A)~=~H_{2r-n}(\Z_2;A,(-1)^{n-r})~;~
f \mapsto (f\otimes f)(\psi_{2r-n})\cr
\widehat{v}_r(\theta)~ :~
H^{n-r}(C)~=~H_0(\Hom_A(C,S^{n-r}A))\cr
\hskip50pt \to \widehat{Q}^n(S^{n-r}A)~=~\widehat{H}^r(\Z_2;A)~;~
f \mapsto (f\otimes f)(\theta_{n-2r})~.\endcases$$
\medskip

\noindent \S5. {\bf Algebraic Poincar\'e complexes}
\medskip

An algebraic Poincar\'e complex is a chain complex with Poincar\'e
duality, such as arises from a compact $n$-manifold or a normal map.
\medskip
 
An {\it $n$-dimensional}
$\cases \text{\it symmetric}\cr \text{\it quadratic}\endcases$ 
{\it {\rm (}Poincar\'e{\rm )} complex over $A$} 
$\cases (C,\phi)\cr (C,\psi)\endcases$ is an
$n$-dimen\-sional $A$-module chain complex $C$ together with an $n$-dimensional
$\cases \text{symmetric}\cr \text{quadratic}\endcases$ structure
$\cases \phi\in (W^{\%}C)\cr \psi\in (W_{\%}C)_n\endcases$ 
(such that $\cases \phi_0 \cr (1+T)\psi_0\endcases:C^{n-*}\to  C$ 
is a chain equivalence). 
\medskip

An {\it $(n+1)$-dimensional}
$\cases \text{\it symmetric} \cr \text{\it quadratic} \cr 
\text{\it {\rm (}symmetric,\t quadratic{\rm )}}\endcases$
{\it {\rm(}Poincar\'e{\rm)} pair over $A$}
$$ \bigg(~f~:~C \to D~,~\cases (\delta\phi,\phi)\cr
(\delta\psi,\psi)\cr (\delta\phi,\psi)\endcases~\bigg)$$
consists of an $n$-dimensional $A$-module chain complex $C$, an
$(n+1)$-dimensional $A$-module chain complex $D$, a chain map 
$f:C\to  D$ and a cycle
$$\cases (\delta\phi,\phi) \in
C(f^{\%}:W^{\%}C\to W^{\%}D)_{n+1}~=~
(W^{\%}D)_{n+1}\oplus (W^{\%}C)_n\cr
(\delta\psi,\psi) \in
C(f_{\%}:W_{\%}C\to W_{\%}D)_{n+1}~=~
(W_{\%}D)_{n+1}\oplus (W_{\%}C)_n\cr
(\delta\phi,\psi) \in
C((1+T)f_{\%}:W_{\%}C\to W^{\%}D)_{n+1}~=~
(W^{\%}D)_{n+1}\oplus (W_{\%}C)_n\endcases$$
(such that the $A$-module chain map $D^{n+1-*}\to C(f)$ defined by 
$$\cases 
(\delta\phi,\phi)_0~=~\pmatrix \delta\phi_0\cr \phi_0f^*\endpmatrix
\cr
(1+T)(\delta\psi,\psi)_0~=~
\pmatrix (1+T)\delta\psi_0 \cr (1+T)\psi_0f^* \endpmatrix\cr
(\delta\phi,(1+T)\psi)_0~=~\pmatrix \delta\phi_0 \cr
(1+T)\psi_0f^*\endpmatrix\endcases~:~
D^{n+1-r} \to C(f)_r~=~D_r\oplus C_{r-1}$$
is a chain equivalence). 
The {\it boundary} of the pair is the $n$-dimensional
$\cases \text{symmetric}\cr \text{quadratic}\cr \text{quadratic}\endcases$
{\rm(}Poincar\'e{\rm)} complex $\cases (C,\phi)\cr (C,\psi)\cr (C,\psi).\endcases$
\medskip

A {\it homotopy equivalence} of $n$-dimensional
$\cases \text{symmetric}\cr \text{quadratic}\endcases$ complexes 
$$\cases (f,\chi)~ :~ (C,\phi)\to (C',\phi')\cr
(f,\xi)~ :~ (C,\psi) \to (C',\psi')\endcases$$
is a chain equivalence $f : C\to C'$ 
together with an equivalence of 
$\cases \text{symmetric}\cr \text{quadratic}\endcases$ structures on $C'$ 
$\cases \chi~:~f^{\%}(\phi) \to \phi'\cr
\xi~:~f_{\%}(\psi) \to \psi'~.\endcases$
There is a similar notion of homotopy equivalence for pairs. 
\medskip

An $n$-dimensional $\cases \text{symmetric}\cr \text{quadratic}\endcases$ 
complex $\cases (C,\phi)\cr (C,\psi)\endcases$ is {\it connected} if 
$$\cases H_0(\phi_0:C^{n-*}\to C)~=~0\cr
H_0((1+T)\psi_0:C^{n-*} \to C)~=~0~.\endcases$$

It was shown in Ranicki [15] that there is a natural one-one
correspondence between the homotopy equivalence classes of connected
$n$-dimensional 
$\cases \text{symmetric}\cr \text{quadratic}\endcases$
complexes over $A$ and the homotopy equivalence classes of
$n$-dimensional\newline
$\cases \text{symmetric}\cr \text{quadratic}\endcases$ Poincar\'e pairs 
over $A$. A connected $n$-dimensional 
$\cases \text{symmetric}\cr \text{quadratic}\endcases$ 
complex $\cases (C,\phi)\cr (C,\psi)\endcases$ determines the
$n$-dimensional $\cases \text{symmetric}\cr \text{quadratic}\endcases$ 
Poincar\'e pair
$$(i_C:\partial C \to C^{n-*},\cases (0,\partial \phi)\cr
(0,\partial \psi)\endcases)$$
defined by

$$\eqalign{
&i_C~=~\pmatrix 0 & 1 \endpmatrix~:~
\partial C_r~=~C_{r+1}\oplus C^{n-r}\to C^{n-r}~,\cr
&d_{\partial C}~=~\cases 
 \pmatrix d_C & (-1)^r\phi_0\cr
      0& (-1)^rd^*_C\endpmatrix~:\cr
 \pmatrix d_C & (-1)^r(1+T)\psi_0\cr
      0& (-1)^rd^*_C\endpmatrix~:\endcases\cr
&\hskip25pt 
\partial C_r~=~C_{r+1}\oplus C^{n-r} \to \partial C_{r-1}~=~C_r\oplus C^{n-r+1}~,}$$

$$\eqalign{&\cases
\partial \phi_0~=~
  \pmatrix (-1)^{n-r-1}T\phi_1 & (-1)^{r(n-r-1)}\cr 
      1 & 0\endpmatrix~:\cr
    \partial \psi_0~=~
  \pmatrix 0& 0 \cr 1&0 \endpmatrix~:
  \endcases\cr
&\hskip25pt \partial C^{n-r-1}~=~C^{n-r}\oplus C_{r+1} \to 
      \partial C_r~=~C_{r+1}\oplus C^{n-r}~,\cr
&\cases \partial \phi_s~=~
 \pmatrix (-1)^{n-r+s-1}T\phi_{s+1} & 0\cr 0& 0\endpmatrix~:\cr
   \hskip25pt
   \partial C^{n-r+s-1}~=~C^{n-r+s}\oplus C_{r-s+1} \to 
        \partial C_r~=~C_{r+1}\oplus C^{n-r}~~(s \ge 1)~,\cr
   \partial \psi_s~=~\pmatrix (-1)^{n-r-s}T\psi_{s-1} & 0\cr 0 & 0 \endpmatrix~ :\cr
\hskip25pt 
   \partial C^{n-r-s-1}~=~C^{n-r-s}\oplus C_{r+s+1} \to 
        \partial C_r~=~C_{r+1}\oplus C^{n-r}~~(s \ge 1)~.\endcases}$$
\medskip

The $(n-1)$-dimensional $\cases \text{symmetric}\cr \text{quadratic}\endcases$
Poincar\'e complex 
$$\cases 
\partial (C,\phi)~=~(\partial C,\partial \phi)\cr
\partial (C,\psi)~=~(\partial C,\partial \psi)\endcases$$
is the {\it boundary} of the connected $n$-dimensional
$\cases \text{symmetric}\cr \text{quadratic}\endcases$
complex $\cases (C,\phi)\cr (C,\psi).\endcases$ The
connected complex $\cases (C,\phi)\cr (C,\psi)\endcases$ is a
Poincar\'e complex if and only if the boundary
$\cases \partial (C,\phi)\cr \partial (C,\psi)\endcases$ is contractible
(= homotopy equivalent to 0). A Poincar\'e complex
$\cases (C,\phi)\cr (C,\psi)\endcases$ is the boundary of an
$(n+1)$-dimensional $\cases \text{symmetric}\cr \text{quadratic}\endcases$ Poincar\'e
pair
$\cases (f:C\to D,(\delta\phi,\phi))\cr (f:C\to D,(\delta\psi,\psi))\endcases$
if and only if it is homotopy equivalent to the boundary of a connected
$(n+1)$-dimensional $\cases \text{symmetric}\cr \text{quadratic}\endcases$ complex. 
\medskip

The $n$-dimensional $\cases \text{symmetric}\cr \text{quadratic}\endcases$ 
Poincar\'e complexes 
$\cases (C,\phi)\cr (C,\psi)\endcases$, $\cases (C',\phi')\cr (C',\psi')\endcases$
are {\it cobordant} if
$\cases (C,\phi)\oplus (C',-\phi')\cr (C,\psi)\oplus (C',-\psi')\endcases$
is the boundary of an $(n+1)$-dimensional\newline
$\cases \text{symmetric}\cr \text{quadratic}\endcases$ Poincar\'e pair
$\cases ((f~f'):C\oplus C'\to D,(\delta\phi,\phi\oplus -\phi'))\cr 
((f~f'):C\oplus C'\to D,(\delta\psi,\psi\oplus -\psi'))~.\endcases$
Homotopy equivalent Poincar\'e complexes are cobordant. 
\medskip

The $\cases \text{\it symmetric}\cr \text{\it quadratic}\endcases$ 
{\it $L$-groups} 
$\cases L^n(A)\cr L_n(A)\endcases$ ($n \ge 0$) are the
cobordism groups of $n$-dimensional $\cases \text{symmetric}\cr \text{quadratic}\endcases$ 
Poincar\'e complexes over $A$. The quadratic $L$-groups
$L_*(A)$ are 4-periodic, with isomorphisms 
$$L_n(A) \to L_{n+4}(A)~ ;~(C,\psi) \mapsto (S^2C,\psi)~~(n \ge 0)~,$$ 
and are just the surgery obstruction groups of Wall [18].  The
symmetric $L$-groups $L^*(A)$ were introduced by Mishchenko [13].  The
corresponding maps in the symmetric $L$-groups
$$L^n(A) \to L^{n+4}(A)~ ;~ (C,\phi) \mapsto (S^2C,\phi)~~(n \ge 0)$$ 
are not isomorphisms in general.  The symmetric and quadratic
$L$-groups are related by an exact sequence
$$\dots \to L_n(A)~\raise4pt\hbox{$1+T \atop \to$}~
L^n(A) \to \widehat{L}^n(A) \to L_{n-1}(A) \to \dots $$
with 
$$1+T ~:~ L_n(A) \to L^n(A)~ ;~ (C,\psi) \mapsto (C,(1+T)\psi)$$
and $\widehat{L}^n(A)$ the relative cobordism group of $n$-dimensional
(symmetric, quadratic) Poinc\-ar\'e pairs over $A$. The relative $L$-groups
$\widehat{L}^*(A)$ are 8-torsion, so that the symmetrization maps
$1+T:L_n(A)\to L^n(A)$ are isomorphisms modulo 8-torsion. If
$\widehat{H}^*(\Z_2;A)=0$ (e.g. if $1/2\in A$) then
$\widehat{L}^*(A)=0$ and the symmetrization maps are isomorphisms.
\medskip

The {\it symmetric construction} of Ranicki [15] is the natural chain map
$$\phi_X~=~1\otimes \Delta~:~C(X) \to W^{\%}C(\widetilde{X})~=~
\Hom_{\Z[\Z_2]}(W,C(\widetilde{X})\otimes_{\Z[\pi]}C(\widetilde{X}))$$
induced by an Alexander-Whitney-Steenrod diagonal chain approximation
$\Delta$, for any space $X$ and any regular cover $\widetilde{X}$, with
$\pi$ the group of covering translations. For $\widetilde{X}=X$ the
mod 2 reduction of the composite
$$H_n(X) \xymatrix{\ar[r]^{\di{\phi_X}}&} Q^n(C(X)) 
\xymatrix{\ar[r]^{\di{v_r}}&} \Hom_{\Z}(H^{n-r}(X),Q^n(S^{n-r}\Z))$$
is given by the $r$th Steenrod square
$$v_r(\phi_X(x))(y)~=~\langle Sq^r(y),x \rangle \in \Z_2~.$$
\medskip

The symmetric signature of Mishchenko [13] is defined for any
$n$-dimensional geometric Poincar\'e complex $X$ to be the symmetric
Poincar\'e cobordism class
$$\sigma^*(X)~=~(C(\widetilde{X}),\phi_X([X])) \in L^n(\Z[\pi_1(X)])~.$$
The symmetric $L$-groups of $\Z$ are given by
$$L^n(\Z)~=~\cases \Z~\text{(signature)} &\text{if $n \equiv 0(\bmod\,4)$}\cr
\Z_2~\text{(deRham invariant)}&\text{if $n\equiv 1(\bmod\,4)$}\cr
0&\text{if $n\equiv 2(\bmod\,4)$}\cr
0&\text{if $n\equiv 3(\bmod\,4)$~.}\endcases$$
\medskip

The {\it quadratic construction} of [15] associates to any 
stable $\pi$-equivariant map $F:\Sigma^{\infty}\widetilde{X}_+ \to 
\Sigma^{\infty}\widetilde{Y}_+$ a natural chain map
$$\psi_F~:~C(X) \to W_{\%}C(\widetilde{Y})~=~
W\otimes_{\Z[\Z_2]}(C(\widetilde{Y})\otimes_{\Z[\pi]}C(\widetilde{Y}))$$
such that 
$$(1+T)\psi_F~=~F^{\%}\phi_X - \phi_YF_*~:~
C(X) \to W^{\%}C(\widetilde{Y})$$
with $\widetilde{X}$ a regular cover of $X$ with group of covering
translations $\pi$, $\widetilde{X}_+=\widetilde{X} \cup \{\text{pt.}\}$,
and similarly for $Y$. For $\widetilde{X}=X$, $\widetilde{Y}=Y$, $\pi=\{1\}$
the mod 2 reduction of the composite
$$H_n(X) \xymatrix{\ar[r]^{\di{\psi_F}}&} Q_n(C(Y)) 
\xymatrix{\ar[r]^{\di{v^r}}&} \Hom_{\Z}(H^{n-r}(Y),Q_n(S^{n-r}\Z))$$
is given by the $(r+1)$th functional Steenrod square
$$v^r(\psi_F(x))(y)~=~\langle Sq_{(\Sigma^{\infty}y)F}^{r+1}
(\Sigma^{\infty}\iota),\Sigma^{\infty}x \rangle \in \Z_2$$
with $\iota \in H^{n-r}(K(\Z_2,n-r);\Z_2)=\Z_2$ the generator.

\medskip

The Wall [18] surgery obstruction of an $n$-dimensional normal map 
$(f,b):M \to X$ was expressed in [15] as the quadratic Poincar\'e cobordism class 
$$\sigma_*(f,b)~=~(C(f^!),e_{\%}\psi_F([X])) \in L_n(\Z[\pi_1(X)])$$
with $F:\Sigma^{\infty} \widetilde{X}_+ \to \Sigma^{\infty}\widetilde{M}_+$
a $\pi_1(X)$-equivariant $S$-dual of
$T(\widetilde{b}):T(\nu_{\widetilde{M}})\to T(\nu_{\widetilde{X}})$ 
inducing the Umkehr chain map
$$f^!~:~ C(\widetilde{X})~\simeq~C(\widetilde{X})^{n-*}
\xymatrix{\ar[r]^{\di{f^*}}&} C(\widetilde{M})^{n-*}~\simeq~C(\widetilde{M})$$
and $e:C(\widetilde{M}) \to C(f^!)$ the inclusion in the algebraic mapping cone.
The quadratic $L$-groups of $\Z$ are given by
$$L_n(\Z)~=~\cases \Z~\text{(signature/8)} &\text{if $n \equiv 0(\bmod\,4)$}\cr
0&\text{if $n\equiv 1(\bmod\,4)$}\cr
\Z_2~\text{(Arf invariant)}&\text{if $n\equiv 2(\bmod\,4)$}\cr
0&\text{if $n\equiv 3(\bmod\,4)$~.}\endcases$$
\medskip

\noindent \S6. {\bf Chain bundles} 

\medskip

A {\it bundle} over a finite-dimensional $A$-module chain complex C is a
0-dimensional hyperquadratic structure on $C^{0-*}$, 
that is a cycle
$$\gamma \in (\widehat{W}^{\% }C^{0-*})_0~,$$
as defined by a collection of $A$-module morphisms 
$$\{\gamma_s\in \Hom_A(C_{r-s},C^{-r})\,\vert\, r,s\in \Z\}$$
such that
$$(-1)^{r+1}d^*_C\gamma_s +(-1)^s\gamma_sd_C +
(-1)^{s-1}(\gamma_{s-1}+(-1)^sT\gamma_{s-1})~=~0~:~C_{r-s+1}\to C^{-r}~.$$ 

An {\it equivalence} of bundles over $C$ 
$$\chi~ :~ \gamma \to \gamma'$$ 
is an equivalence of hyperquadratic structures, 
as defined by a collection of $A$-module morphisms 
$$\{\chi_s\in \Hom_A(C_{r-s-1},C^{-r})\,\vert\, r,s\in \Z\}$$
such that
$$\gamma'_s -\gamma_s~=~(-1)^{r+1}d^*_C\chi_s +
(-1)^s\chi_sd_C+(-1)^s(\chi_{s-1}+(-1)^sT\chi_{s-1})~ :~C_{r-s}\to C^{-r}~.$$
Thus
$$\widehat{Q}^0(C^{0-*})~=~H_0(\widehat{W}^{\%}C^{0-*})$$
is the abelian group of equivalence classes of bundles
over $C$. 
\medskip

A {\it chain bundle} over $A$ $(C,\gamma)$ is a
finite-dimensional $A$-module chain complex $C$ together with a bundle
$\gamma\in (\widehat{W}^{\% }C^{0-*})_0$. 
\medskip

Given a chain bundle $(C,\gamma)$ over $A$ and an $A$-module chain map 
$f:B\to C$ define the {\it pullback} chain bundle $(B,f^*\gamma)$ using the
image of $\gamma$ under the $\Z$-module chain map 
$$\widehat{f}^*~:~\widehat{W}^{\% }C^{0-*} \to \widehat{W}^{\%}B^{0-*}$$
induced by the dual $A$-module chain map $f^*:C^{0-*}\to B^{0-*}$.
The equivalence class of the pullback
bundle $f^*\gamma$ depends only on the chain homotopy class of the
chain map $f$, by the chain homotopy invariance of the $Q$-groups. 
\medskip

A {\it map} of chain bundles over $A$ 
$$(f,\chi)~ :~ (C,\gamma) \to (C',\gamma')$$
is a chain map $f:C\to C'$ together with an equivalence of bundles over $C$
$$\chi~:~\gamma \to f^*\gamma'~.$$ 
The {\it composite} of chain bundle maps 
$$(f,\chi)~ :~(C,\gamma) \to (C',\gamma')~,~
(f',\chi')~ :~(C',\gamma') \to (C'',\gamma'')$$
is the chain bundle map
$$(f',\chi')(f,\chi)~=~(f'f,\chi+\widehat{f^*}^{\%}(\chi'))~:~
(C,\gamma) \to (C'',\gamma'')~.$$
A {\it homotopy of chain bundle maps}
$$(g,\eta)~:~(f,\chi)~\simeq ~(f',\chi')~:~(C,\gamma) \to (C',\gamma')$$
is a chain homotopy 
$$g~ :~ f~ \simeq~ f'~:~ C \to C'$$
together with an equivalence of 1-dimensional hyperquadratic
structures on $C^{0-*}$ 
$$\eta~ :~\chi - \chi' +(g^*;f^*,f^{\prime *})^{\% }(\gamma') \to 0~ .$$
A map of chain bundles $(f,\chi):(C,\gamma)\to (C',\gamma')$ is an {\it
equivalence} if there exists a homotopy inverse. This happens
precisely when $f:C\to C'$ is a chain equivalence, in which case any
chain homotopy inverse
$$f'~=~f^{-1}~ :~ C' \to C$$ 
can be used to define a homotopy inverse
$$(f',\chi')~ :~ (C',\gamma') \to (C,\gamma)~.$$
Given a chain bundle $(B,\beta)$ over $A$ and a finite-dimensional 
$A$-module chain complex $C$
use the pullback construction to define abelian group morphisms
$$I_{\beta}~:~H_n(B^t\otimes _AC) \to \widehat{Q}^n(C)~ ;~
f \mapsto \widehat{f}^{\%}(S^n\beta)~,$$ 
using the slant isomorphism 
$$B^t\otimes _AC \to  \Hom_A(B^{-*},C)~;~
x\otimes y \mapsto (f \mapsto \overline{f(x)}\,.\,y)$$ 
to identify a cycle $f\in (B^t\otimes _AC)_n$ with a chain map
$f:B^{n-*}\to C$. Weiss [19] developed an algebraic analogue of the
representation theorem of Brown [3] to obtain for any ring with
involution $A$ the existence of a directed system
$\{(B(r),\beta(r))\,\vert\, r\ge 0\}$ of chain bundles over $A$ and
chain bundle maps 
$$(B(r),\beta(r)) \to (B(r+1),\beta(r+1))$$
such that the abelian group morphisms 
$$\varinjlim\limits_r\, I_{\beta_r}~ :~\varinjlim\limits_r\,
H_n(B(r)^t\otimes_AC) \to \widehat{Q}^n(C)$$
are isomorphisms for any finite-dimensional $A$-module chain complex
$C$.  In general, the direct limit $A$-module chain complex
$$B(\infty)~=~\varinjlim\limits_r\, B(r)$$ 
is not finite-dimensional. As in [19] we shall ignore
this inconvenience and treat $B(\infty)$ as if it were finite-dimensional,
so that there is defined the {\it universal chain bundle} over $A$
$$(B(\infty),\beta(\infty))~=~\varinjlim\limits_r\, (B(r),\beta(r))~,$$ 
with the universal property that for any finite-dimensional
$A$-module chain complex $C$ the abelian group morphisms 
$$I_{\beta(\infty)}~:~H_n(B(\infty)^t\otimes _AC) \to 
\widehat{Q}^n(C)$$
are isomorphisms. In particular, there are defined isomorphisms 
$$H_0(\Hom_A(C,B(\infty))) \to \widehat{Q}^0(C^{0-*})~;~
f \mapsto f^*(\beta(\infty))$$
for any finite-dimensional $C$. Thus every chain bundle $(C,\gamma)$ 
has a classifying map 
$$(f,\chi)~:~(C,\gamma) \to (B(\infty),\beta(\infty))$$
and the equivalence classes of bundles 
$\gamma\in (\widehat{W}^{\% }C^{0-*})_0$ over $C$ are in
one-one correspondence with the chain homotopy classes of chain maps
$f:C\to B(\infty)$. 
\medskip

The {\it Wu classes} of a chain bundle
$(C,\gamma)$ are the Wu classes of $\gamma$, the $A$-module morphisms
$$\widehat{v}_r(\gamma)~:~H_r(C) \to \widehat{H}^r(\Z_2;A)~ ;~
x \mapsto \langle \gamma_{-2r},x\otimes x\rangle~.$$
The universal chain bundle $(B(\infty),\beta(\infty))$
is characterized by the property that the Wu
classes define $A$-module isomorphisms 
$$\widehat{v}_r(\gamma)~:~H_r(C) \to \widehat{H}^r(\Z_2;A)~~(r \ge 0)~.$$

For example, if $A=\Z$ the chain bundle $(B(\infty),\beta(\infty))$
defined by 
$$\eqalign{
&d_{B(\infty)}~=~\cases 2&\text{if $r$ is odd} \cr 0&\text{if $r$ is even}
\endcases ~:~B(\infty)_r~=~\Z \to B(\infty)_{r-1}~=~\Z~,\cr
&\beta(\infty)_s~=~\cases 1&\text{if $2r=s$} \cr 0&\text{otherwise}\endcases~:~
B(\infty)_{r-s}~=~\Z \to B(\infty)^{-r}~=~\Z}$$
is universal, with the Wu classes defining isomorphisms
$$\widehat{v}_r(\beta(\infty))~:~H_r(B(\infty))
\to \widehat{H}^r(\Z_2;\Z)~=~\cases \Z_2&\text{if $r$ is even} \cr
0&\text{if $r$ is odd~.}\endcases$$
\medskip

The symmetric and quadratic constructions of Ranicki [15] were extended
in Weiss [19] and Ranicki [16,\t 2.5]\t : a spherical fibration $\nu:X
\to BG(k)$ determines a chain bundle $(C(\widetilde{X}),\gamma)$ over
$\Z[\pi_1(X)]$, and there is defined a natural transformation of exact
sequences from the certain exact sequence of Whitehead [20]

$$\xymatrix@C-8pt{\dots \ar[r] 
& \Gamma_{n+k+1}(T(\nu)) \ar[r] \ar[d]
& \pi_{n+k}(T(\nu)) \ar[r]^-{\di{h}} \ar[d]
& \dot H_{n+k}(T(\nu)) \ar[r] \ar[d]^-{\di{\phi_XU}}
& \Gamma_{n+k}(T(\nu)) \ar[r] \ar[d]&\dots \\
\dots \ar[r]
& \hbox{$\widehat{Q}^{n+1}(C(\widetilde X))$} \ar[r] 
& Q_n(C(\widetilde X),\gamma) \ar[r] 
& Q^n(C(\widetilde X)) \ar[r]^-{\di{J_{\gamma}}}
& \hbox{$\widehat{Q}^n(C(\widetilde X))$} \ar[r]&\dots}$$

\medskip
\noindent with $h$ the Hurewicz map and $U:\dot H_{n+k}(T(\nu)) \to H_n(X)$
the Thom isomorphism. The topological Wu classes of $\nu$ 
are the algebraic Wu classes of the induced chain bundle 
$(C(X;\Z_2),1\otimes \gamma)$ over $\Z_2$
$$v_*(\nu)~=~\widehat{v}_*(1\otimes \gamma) \in 
H^*(X;\Z_2)~=~\text{Hom}_{\Z_2}(H_*(X;\Z_2),\Z_2)~.$$

\noindent \S7. {\bf Normal complexes}
\medskip

An $n$-dimensional normal space $(X,\nu_X,\rho_X)$ (Quinn [14]) is a
finite $n$-dimensional $CW$ complex $X$ together with a $(k-1)$-spherical
fibration $\nu_X:X\to BG(k)$ and a map
$\rho_X:S^{n+k}\to T(\nu_X)$ to the Thom space of $\nu_X$. An
$n$-dimensional geometric Poincar\'e complex $X$ has a unique equivalence
class of normal structures $(\nu_X,\rho_X)$, with $\nu_X$ the
Spivak normal fibration and $\rho_X$ representing the fundamental
class $[X]\in H_n(X)$. A normal complex is the algebraic analogue of
a normal space, consisting of a symmetric complex with normal chain
bundle. 
\medskip

A {\it normal structure} $(\gamma,\theta)$ on an
$n$-dimensional symmetric complex $(C,\phi)$ is a bundle 
$\gamma\in (\widehat{W}^{\% }C^{0-*})_0$ together with an equivalence of
$n$-dimensional hyperquadratic structures on $C$
$$ \theta~ :~J(\phi) \to (\widehat\phi_0)^{\% }(S^n\gamma)~,$$
as defined by a chain $\theta\in (\widehat{W}^{\% }C)_{n+1}$ 
such that 
$$J(\phi)- (\widehat\phi_0)^{\% }(S^n\gamma)~=~d\theta \in 
(\widehat{W}^{\% }C)_n~.$$ 
The Wu classes of $\phi$ and $\gamma$ are then related by a commutative diagram 
$$\xymatrix@C+10pt{
H^{n-r}(C) \ar[r]^-{\di{\phi_0}} \ar[d]_-{\di{v_r(\phi)}}
& H_r(C) \ar[d]_-{\hbox{$\widehat{v}_r(\gamma)$}} \\
H^{n-2r}(\Z_2;A,(-1)^{n-r}) \ar[r]^-{\di{J}} & \hbox{$\widehat{H}^r(Z_2;A)$}~.
}$$

An {\it equivalence}
of $n$-dimensional normal structures on $(C,\phi)$ 
$$(\chi,\eta)~ :~(\gamma,\theta) \to (\gamma',\theta')$$
is an equivalence of bundles $\chi:\gamma\to \gamma'$ 
together with an equivalence of $(n+1)$-dimensional hyperquadratic 
structures on $C$
$$\eta~ :~ \theta -\theta' + (\widehat\phi_0)^{\%}(S^n\chi) \to 0 ~.$$

An $n$-dimensional symmetric Poincar\'e complex $(C,\phi\in (W^{\%}C)_n)$ 
has a unique equivalence class of normal structures
$(\gamma,\theta)$, with the equivalence class of bundles
$[\gamma]\in \widehat{Q}^0(C^{0-*})$ the image of the
equivalence class of symmetric structures $[\phi]\in Q^n(C)$ under
the composite 
$$Q^n(C)~\raise4pt\hbox{$J \atop \to$}~\widehat{Q}^n(C)~
\raise4pt\hbox{$((\phi_0)^{\%})^{-1} \atop \to$}~\widehat{Q}^n(C^{n-*})~
\raise4pt\hbox{$(S^n)^{-1} \atop \to$}~\widehat{Q}^0(C^{0-*})~.$$ 
If $(\gamma,\theta)$, $(\gamma',\theta')$ are two such normal
structures on $(C,\phi)$ there exists an equivalence of bundles
$\chi:\gamma\to \gamma'$. As $\phi_0:C^{n-*}\to C$ is a chain
equivalence the cycle 
$$\theta - \theta' + (\phi_0)^{\%}(S^n\chi) \in (\widehat{W}^{\% }C)_{n+1}$$
is such that there exist a cycle 
$\lambda\in \widehat{W}^{\%}C^{n-*})_{n+1}$ 
and a chain $\mu\in (\widehat{W}^{\%}C)_{n+2}$ such that 
$$\theta - \theta' + (\phi_0)^{\%}(S^n\chi)~=~
(\phi_0)^{\% }(\lambda) + d\mu \in (\widehat{W}^{\% }C)_{n+1}~.$$
There is now defined an equivalence of normal structures on $(C,\phi)$
$$(\chi-(S^n)^{-1}(\lambda),\mu)~ :~
(\gamma,\theta) \to (\gamma',\theta')~ .$$

An {\it $n$-dimensional normal {\rm(}Poincar\'e{\rm)} complex} over $A$
$(C,\phi,\gamma,\theta)$ is an $n$-dimensio\-nal symmetric (Poincar\'e)
complex $(C,\phi)$ together with a normal structure
$(\gamma,\theta)$. Symmetric Poincar\'e complexes are regarded as
normal Poincar\'e complexes by choosing a normal structure in the
unique equivalence class. 
\medskip

An $n$-dimensional normal complex
$(C,\phi,\gamma,\theta)$ is {\it connected} if the $n$-dimensional
symmetric complex $(C,\phi)$ is connected, that is
$$H_0(\phi_0:C^{n-*}\to C)~=~0~.$$

The correspondence described in \S5 between the homotopy equivalence 
classes of connected $n$-dimensional 
$\cases \text{symmetric}\cr \text{quadratic}\endcases$ complexes and
those of $n$-dimensional $\cases \text{symmetric}\cr \text{quadratic}\endcases$ 
Poincar\'e pairs has the following generalization to connected normal
complexes and (symmetric, quadratic) Poincar\'e pairs. 
\medskip

A connected $n$-dimensional normal complex $(C,\phi,\gamma,\theta)$
determines the $n$-dimen\-sional (symmetric, quadratic) Poincar\'e pair
$$(i_C: \partial C \to C^{n-*},(\delta \phi,\psi))$$
defined by
$$\eqalign{
&i_C~=~(0~1)~: ~\partial C_r~=~C_{r+1} \oplus C^{n-r}\to C^{n-r}~,\cr
&d_{\partial C}~=~\pmatrix d_C & (-1)^r \phi_0 \cr
0 & (-1)^rd^*_C \endpmatrix~: \cr
& \indent \partial C_r~=~C_{r+1} \oplus C^{n-r}
\to\partial C_{r-1}~=~C_r \oplus C^{n-r+1}~,\cr
&\psi_0~=~\pmatrix \chi_0 & 0 \cr
1+ \gamma_{-n}\phi^*_0 & \gamma^*_{-n-1}\endpmatrix~: \cr
&\indent \partial C^r~=~C^{r+1} \oplus C_{n-r}
\to\partial C_{n-r-1}~=~C_{n-r} \oplus C^{r+1}~,\cr
&\psi_s~=~\pmatrix \chi_{-s} & 0 \cr
\gamma_{-n-s}\phi^*_0 & \gamma^*_{-n-s-1}\endpmatrix~: \cr
&\indent \partial C^r~=~C^{r+1} \oplus C_{n-r}
\to\partial C_{n-r-s-1}~=~C_{n-r-s} \oplus C^{r+s+1}~(s \ge 1)~,\cr
&\delta \phi_s~=~\gamma_{-n-s}~: ~C_r \to C^{n-r+s}~(s\ge 
0)~.}$$
The $(n-1)$-dimensional quadratic Poincar\'e complex
$$\partial (C,\phi,\gamma,\theta)~=~(\partial C,\psi)$$ is the
{\it quadratic boundary} of the connected $n$-dimensional normal complex
$(C,\phi,\gamma,\theta)$. (Compare with the definition in \S6 of
the boundary $(n-1)$-dimensional
$\cases \text{symmetric}\cr \text{quadratic}\endcases$ Poincar\'e complex
$\cases \partial (C,\phi)\cr \partial (C,\psi)\endcases$ of a connected
$n$-dimensional $\cases \text{symmetric}\cr \text{quadratic}\endcases$ complex
$\cases (C,\phi)\cr (C,\psi)\endcases$). Conversely, given an
$n$-dimensional (symmetric, quad\-ratic) Poincar\'e pair 
$(f:C\to D,(\delta\phi,\psi))$ there is defined a
connected $n$-dimensional normal complex $(C(f),\phi,\gamma,\theta)$
with the symmetric structure
$$\eqalign{\phi_s~=~& \cases  
\pmatrix \delta \phi_0 & 0 \cr\noalign{\vskip2pt} 
(1+T)\psi_0f^* & 0\endpmatrix &\text{if $s=0$} \cr
\pmatrix \delta \phi_1 & 0 \cr\noalign{\vskip2pt} 
0 & (1+T) \psi_0\endpmatrix &\text{if $s=1$} \cr
\pmatrix \delta \phi_s & 0 \cr 0 & 0\endpmatrix &\text{if $s \ge 2$}\endcases\cr 
&: ~C(f)^r~=~D^r \oplus C^{r-1} \to C(f)_{n-r+s}~=~D_{n-r+s} \oplus 
C_{n-r+s-1}~.}$$
The normal structure $(\gamma,\chi)$ is determined up to equivalence by the 
Poincar\'e duality, with $\gamma \in \widehat Q^0(D^{-*})$ the image of
$(\delta \phi /(1+T)\psi) \in Q^n(C(f))$ under the composite
$$Q^n(C(f)) \xymatrix@C+60pt{\ar[r]^-{\di{((\delta \phi_0,(1+T)\psi_0)^{\%})^{-1}}}&} 
Q^n(D^{n-*}) \xymatrix{\ar[r]^{\di{J}}&} \widehat Q^n(D^{n-*})~
\xymatrix@C+10pt{\ar[r]^{\di{S^{-n}}}&} \widehat Q^0(D^{-*})~.$$
The composite isomorphism
$$\widehat{Q}^0(C(f)^{0-*}) \xymatrix{\ar[r]^{\di{S^n}}&}
\widehat{Q}^n(C(f)^{n-*}) 
\xymatrix@C+40pt{\ar[r]^-{\di{(\delta \phi_0,(1+T)\psi_0)^{\%}}}&} 
\widehat{Q}^n(D)$$
sends the equivalence class $[\gamma]\in \widehat{Q}^0(C(f)^{0-*})$
to the element $\alpha\in \widehat{Q}^n(D)$ 
represented by
$$\alpha_s~=~\cases \delta\phi_s&\text{if $s\ge 0$}\cr
f\psi_{-s-1}f^*&\text{if $s \le -1$}\endcases~:~D^r \to D_{n-r+s}~.$$
There is thus a natural one-one
correspondence between the homotopy equivalence classes of connected
$n$-dimensional normal complexes over $A$ and the homotopy equivalence
classes of $n$-dimensional (symmetric, quadratic) Poincar\'e pairs over $A$. 
In \S8 below this correspondence will be used to identify the
cobordism group $\widehat{L}^n(A)$ of $n$-dimensional
(symmetric, quadratic) Poincar\'e pairs over $A$ with the cobordism group
of $n$-dimensional normal complexes over $A$. 
\medskip

Let $(B,\beta)$ be a chain bundle over $A$.
A {\it normal $(B,\beta)$-structure}
$(\gamma,\theta,f,\chi)$ on an $n$-dimensional symmetric complex
$(C,\phi)$ over $A$ is a normal structure $(\gamma,\theta)$ on
$(C,\phi)$ together with a chain bundle map 
$$(f,\chi)~ :~(C,\gamma) \to (B,\beta)~.$$
There are also the corresponding relative notions of {\it normal} 
$(B,\beta)$-structure on symmetric and
(symmetric, quadratic) pairs. For the universal chain bundle
$(B(\infty),\allowbreak\beta(\infty))$ over $A$ a normal
$(B(\infty),\beta(\infty))$-structure $(\gamma,\theta,f,\chi)$ on a
symmetric complex $(C,\phi)$ is to all intents and purposes the same as
a normal structure $(\gamma,\theta)$. 
\medskip

A normal (0,0)-structure $(\gamma,\theta,0,\chi)$ on an $n$-dimensional
symmetric complex $(C,\phi)$ determines an equivalence to 0 of the
hyperquadratic structure $J(\phi)\in (\widehat{W}^{\%}C)_n$
$$\xi~=~\theta + \phi_0^{\%}(S^n\chi)~:~J(\phi) \to 0~ .$$
Such an equivalence $\xi:J(\phi)\to 0$ consists of a quadratic
structure $\psi\in (W_{\%}C)_n$ and an equivalence of symmetric
structures 
$$\eta~:~(1+T)\psi \to \phi~ ,$$
with
$$\eqalign{
&\psi_s~=~\xi_{-s-1} \in \Hom_A(C^*,C)_{n-s}~~(s \ge 0)~,\cr
&\eta_s~=~\xi_s \in \Hom_A(C^*,C)_{n+s+1}~~(s\ge 0)~.}$$
Thus a normal (0,0)-structure on a symmetric complex $(C,\phi)$ 
is to all intents and purposes an equivalence of the symmetric structure 
$\phi$ to $(1+T)\psi$ for some quadratic structure $\psi$ on $C$. 
\medskip

An {\it $n$-dimensional $(B,\beta)$-normal {\rm(}Poincar\'e{\rm)} complex}
$(C,\phi,\gamma,\theta,f,\chi)$ is an $n$-dimen\-sion\-al symmetric
(Poincar\'e) complex $(C,\phi)$ together with a normal
$(B,\beta)$-structure $(\gamma,\theta,\phi,\chi)$. 
\medskip

In \S8 below the cobordism group $\widehat{L}\langle B,\beta\rangle^n(A)$ of
$n$-dimensional $(B,\beta)$-normal complexes over $A$ will be identified
with the twisted quadratic group $Q_n(B,\beta)$ (introduced by Weiss
[19]) of equivalence classes of pairs $(\phi,\theta)$ such that
$(B,\phi,\beta,\theta,1,0)$ is an $n$-dimensional $(B,\beta)$-normal
complex. 
\medskip

An {\it $n$-dimensional symmetric structure}
$(\phi,\theta)$ on a chain bundle $(C,\gamma)$ is an $n$-dimen\-sional 
symmetric structure $\phi\in (W^{\% }C)_n$ together with an equivalence 
of $n$-dimensional hyperquadratic structures on $C$ 
$$\theta~:~J(\phi)\to (\phi_0)^{\% }(S^n\gamma)~,$$
as defined by a chain 
$\theta\in (\widehat{W}^{\% }C)_{n+1}$ such that 
$$J(\phi) - (\phi_0)^{\% }(S^n\gamma)~=~d(\theta) \in 
(\widehat{W}^{\% }C)_n ~.$$
Thus $(C,\phi)$ is an $n$-dimensional
symmetric complex with normal structure $(\gamma,\theta)$. 
\medskip

An {\it equivalence} of $n$-dimensional symmetric structures on
$(C,\gamma)$ 
$$(\xi,\eta)~ :~ (\phi,\theta) \to (\phi',\theta')$$
is defined by an equivalence of symmetric structures 
$\xi : \phi\to \phi'$ together with an equivalence of hyperquadratic
structures on $C$ 
$$\eta~ :~ \theta - \theta' + J(\xi) +
(\xi_0;\phi_0,\phi'_0)^{\%}(S^n\gamma) \to 0~ ,$$
as defined by chains $\xi\in (W^{\% }C)_{n+1}$, 
$\eta\in (\widehat{W}^{\% }C)_{n+2}$ such that 
$$\eqalign{
&\phi' - \phi~=~d(\xi) \in (W^{\% }C)_nC^{-*}\cr
&\theta' - \theta + J(\xi) +
(\xi_0;\phi_0,\phi'_0)^{\%}(S^n\gamma)~=~d(\eta) \in 
(\widehat{W}^{\% }C)_{n+1}~.}$$
The {\it twisted quadratic $Q$-group} $Q_n(C,\gamma)$ is the
abelian group of equivalence classes of $n$-dimensional symmetric
structures on a chain bundle $(C,\gamma)$, with addition by
$$(\phi,\theta) +
(\phi',\theta')~=~(\phi+\phi',\theta+\theta'+[\phi_0,\phi'_0](S^n\gamma))
\in Q_n(C,\gamma) ~.$$
The twisted quadratic $Q$-groups
$Q_*(C,\gamma)$ fit into an exact sequence of abelian groups
$$\dots \to \widehat{Q}^{n+1}(C)~
\raise4pt\hbox{$\displaystyle{H_{\gamma}}\atop \to$}~Q_n(C,\gamma)~
\raise4pt\hbox{$\displaystyle{N_{\gamma}} \atop \to$}~Q^n(C)~
\raise4pt\hbox{$\displaystyle{J_{\gamma}} \atop \to$}~
\widehat{Q}^n(C) \to \dots$$
with the morphisms 
$$\eqalign{&H_{\gamma}~:~\widehat{Q}^{n+1}(C) \to Q_n(C,\gamma)~ ;~ 
\theta \mapsto (0,\theta) ~,\cr
&N_{\gamma}~ :~Q_n(C,\gamma) \to Q^n(C)~ ;~
(\phi,\theta) \mapsto \phi~,\cr
&J_{\gamma}~:~ Q^n(C) \to \widehat{Q}^n(C)~;~
\phi \mapsto J(\phi)-(\phi_0)^{\%}(S^n\gamma)}$$
induced by simplicial maps. In the untwisted case
$\gamma=0$ there is defined an isomorphism of exact sequences
$$\xymatrix@C+10pt{
\dots \ar[r] & 
\hbox{$\widehat{Q}^{n+1}(C)$} \ar@{=}[d] \ar[r]^-{\di{H}} & 
Q_n(C) \ar[d]^-{\di{\cong}} \ar[r]^-{\di{1+T}} &
Q^n(C) \ar@{=}[d] \ar[r]^-{\di{J}} &
\hbox{$\widehat{Q}^n(C)$} \ar@{=}[d] \ar[r] & \dots\\
\dots \ar[r] & 
\hbox{$\widehat{Q}^{n+1}(C)$}  \ar[r]^-{\di{H_0}} & 
Q_n(C,0) \ar[r]^-{\di{N_0}} &
Q^n(C) \ar[r]^-{\di{J_0}} &
\hbox{$\widehat{Q}^n(C)$} \ar[r] & \dots}$$
with 
$$Q_n(C) \to Q_n(C,0)~;~\psi \mapsto ((1+T)\psi,\theta)~ ,~
\theta_s~=~\cases \psi_{-s-1}&\text{if $s \le -1$} \cr
0&\text{if $s \ge 0$}.\endcases$$
The twisted quadratic groups $Q_*(C,\gamma)$ are covariant in $(C,\gamma)$. 
Given a map of chain bundles $(f,\chi) : (C,\gamma)\to (C',\gamma')$ 
and an $n$-dimensional symmetric structure $(\phi,\theta)$ on $(C,\gamma)$ 
define an $n$-dimensional symmetric structure on $(C',\gamma')$ 
$$(f,\chi)_{\%}(\phi,\theta)~=~
(f^{\% }(\phi),\widehat{f}^{\%}(\theta)+(f\phi_0)^{\% }(S^n\chi))~.$$ 
The resulting morphisms of the twisted quadratic $Q$-groups 
$$(f,\chi)_{\% }~ :~Q_n(C,\gamma) \to Q_n(C',\gamma')$$
depend only on the homotopy class of $(f,\chi)$. 
There is defined a morphism of exact sequences
$$\xymatrix{
\dots \ar[r] & 
\hbox{$\widehat{Q}^{n+1}(C)$} 
\ar[d]_-{\di{ \hbox{$\widehat{f}^{\%}$} }} \ar[r]^-{\di{H_{\gamma}}} & 
Q_n(C,\gamma) \ar[d]_-{\di { (f,\chi)_{\%} } } \ar[r]^-{\di{N_{\gamma}}} &
Q^n(C) \ar[d]_-{\di{f^{\%}}} \ar[r]^-{\di{J_{\gamma}}} &
\hbox{$\widehat{Q}^n(C)$} \ar[d]_-{\di{\hbox{$\widehat{f}^{\%}$}}}
\ar[r] & \dots\\
\dots \ar[r] & 
\hbox{$\widehat{Q}^{n+1}(C')$}  \ar[r]^-{\di{H_{\gamma'}}} & 
Q_n(C',\gamma') \ar[r]^-{\di{N_{\gamma'}}} &
Q^n(C') \ar[r]^-{\di{J_{\gamma'}}} &
\hbox{$\widehat{Q}^n(C')$} \ar[r] & \dots}$$
which is an isomorphism if $(f,\chi)$ is an equivalence. 
\medskip

The {\it characteristic element} of an $n$-dimensional $(B,\beta)$-normal 
complex $(C,\phi,\gamma,\allowbreak \theta,f,\chi)$ is defined by 
$$(f,\chi)_{\%}(\phi,\theta) \in Q_n(B,\beta)~.$$ 
In \S8 the cobordism class of a $(B,\beta)$-normal complex will be 
identified with the characteristic element. 
\medskip

A {\it map of $n$-dimensional $\cases\text{\it normal}
\cr \text{\it  $(B,\beta)$-normal}\endcases$ complexes}
$$\cases (f,\xi,\chi,\eta)~:~
(C,\phi,\gamma,\theta) \to (C',\phi',\gamma',\theta')\cr
(f,\xi,\chi,\eta,h,\mu)~:~(C,\phi,\gamma,\theta,g,\lambda) \to 
(C',\phi',\gamma',\theta',g',\lambda')\endcases$$
consists of
\medskip
 
{\parskip=2pt
\parindent=23pt
\item{(i)} a chain map $f:C\to C'$, 
\item{(ii)} an equivalence
$\xi:f^{\% }(\phi)\to \phi'$ of $n$-dimensional symmetric structures
on $C'$, 
\item{(iii)} an equivalence $\chi:\gamma\to f^*\gamma'$ of bundles on $C$, 
\item{(iv)} an equivalence of $(n+1)$-dimensional hyperquadratic structures 
on $C'$}
$$\eta~:~J(\xi) + \theta' -
\widehat{f}^{\%}(\theta) + (\xi_0;f^{\%}\phi_0,\phi'_0)^{\% }(S^n\gamma') 
+(f\phi_0)^{\%}(S^n\chi) \to 0~,$$
\noindent and in the $(B,\beta)$-normal case also 
\medskip

{\parskip=2pt
\parindent=23pt
\item{(v)} a homotopy of bundle maps}
$$(h,\mu)~:~(g,\lambda)~\simeq ~(g',\lambda')(f,\chi)~:~(C,\gamma) \to (B,\beta)~ .$$
Note that $(C,\phi,\gamma,\theta,g,\lambda)$ and
$(C',\phi',\gamma',\theta',g',\lambda')$ have the same
characteristic element 
$$(g,\lambda)_{\%}(\phi,\theta)~=~(g',\lambda')_{\% }(\phi',\theta') \in 
Q_n(B,\beta)~.$$ 
\medskip

It is convenient for computational purposes to describe the behaviour
of the twisted quadratic groups under direct sum.  The direct sum of
chain bundles $(C,\gamma)$, $(C',\gamma')$ is the chain bundle
$$(C,\gamma)\oplus (C',\gamma')~=~(C\oplus C',\gamma\oplus \gamma')~.$$
Let 
$$\eqalign{
&i~=~\pmatrix 1 \cr 0\endpmatrix~:~C \to C\oplus C'~,~
i'~=~\pmatrix 0 \cr 1\endpmatrix~:~C' \to C\oplus C'~,\cr
&j~=~\pmatrix 1 & 0\endpmatrix~:~C\oplus C' \to C~ ,~
j'~=~\pmatrix 0 & 1\endpmatrix~:~C\oplus C' \to C'~ .}$$
The twisted quadratic groups of the direct sum are such that there is
defined a long exact sequence
$$\eqalign{&\dots \to Q_n(C,\gamma)\oplus Q_n(C',\gamma')~
\raise4pt\hbox{$i_* \atop \to$}~Q_n(C\oplus C',\gamma\oplus \gamma')~
\raise4pt\hbox{$j_* \atop \to$}~H_n(C^t\otimes_AC')\cr
&\hskip50pt 
\raise4pt\hbox{$k_* \atop \to$}~
Q_{n-1}(C,\gamma)\oplus Q_{n-1}(C',\gamma')~\raise4pt\hbox{$i_* \atop \to$}~
Q_{n-1}(C\oplus C',\gamma\oplus \gamma') \to \dots}$$
with 
$$\eqalign{
&i_*~=~\pmatrix i_{\%}&i'_{\%}\endpmatrix~:~
Q_n(C,\gamma)\oplus Q_n(C',\gamma') \to Q_n(C\oplus C',\gamma\oplus \gamma')~,
\cr
&j_*~ :~
Q_n(C\oplus C',\gamma\oplus \gamma') \to H_n(C^t\otimes _AC')~;~
(\phi,\theta) \mapsto (j\otimes j')\phi_0~,\cr
&k_*~:~H_n(C^t\otimes_AC') \to Q_{n-1}(C,\gamma)\oplus Q_{n-1}(C',\gamma')~ ;\cr
&\hskip25pt (f:C^{n-*}\to C') \mapsto 
((0,\widehat{f}^{\%}(S^n\gamma')),(0,-\widehat{f^*}^{\%}(S^n\gamma))~ .}$$
For $\gamma=0$ and $\gamma'=0$ 
the long exact sequence collapses into split exact
sequences of the untwisted quadratic $Q$-groups 
$$0 \to Q_n(C)\oplus Q_n(C') \to Q_n(C\oplus C') \to 
H_n(C^t\otimes _AC') \to 0~.$$

\medskip

\noindent \S8. {\bf Normal cobordism}
\medskip

Given a $k$-plane vector bundle $\nu:X\to BO(k)$ over a space $X$ let
$\Omega_n(X,\nu)$ $(n\ge 0)$ denote the bordism groups of bundle
maps 
$$(f,b)~ :~ (M^n,\nu_M) \to (X,\nu)$$ 
with $M^n$ a smooth closed $n$-manifold and $\nu_M:M\to BO(k)$ 
the normal bundle of an embedding $M^n \subset S^{n+k}$ (Lashof [9]). 
The Thom space of $\nu_M$ is given by
$$T(\nu_M)~=~E(\nu_M)/\partial E(\nu_M)$$ 
with $E(\nu_M)$ the tubular neighbourhood of $M^n$ in $S^{n+k}$, 
so that there is defined a collapse map
$$\rho_M~ :~ S^{n+k} \to 
S^{n+k}/(S^{n+k} \backslash E(\nu_M))~=~
E(\nu_M)/\partial E(\nu_M)~=~T(\nu_M)~.$$
The Pontrjagin-Thom isomorphism 
$$\eqalign{&\Omega_n(X,\nu) \to \pi_{n+k}(T(\nu))~;\cr
&(f:M^n\to X,b:\nu_M\to \nu) \mapsto 
(T(b)(\rho_M):S^{n+k}\,\raise4pt\hbox{$\rho_M \atop \to$}\,T(\nu_M)\,
\,\raise4pt\hbox{$T(b) \atop \to$}\,T(\nu))}$$
has inverse 
$$\eqalign{
&\pi_{n+k}(T(\nu)) \to \Omega_n(X,\nu)~;\cr
&(\rho:S^{n+k}\to T(\nu)) \mapsto 
(f=\rho\vert:M^n=\rho^{-1}(X)\to X,b:\nu_M\to \nu)~,}$$
using smooth
transversality to choose a representative $\rho$ transverse regular at
the zero section $X \subset T(\nu)$. 
\medskip

Given a $(k-1)$-spherical
fibration $\nu:X\to BG(k)$ over a space $X$ let $\Omega^N_n(X,\nu)$
(resp. $\Omega^P_n(X,\nu)$) denote the bordism group of
fibration maps 
$$(f,b)~:~(M^n,\nu_M) \to (X,\nu)$$
with $(M^n,\nu_M:M\to BG(k),\rho_M:S^{n+k}\to T(\nu_M))$ an 
$n$-dimensional normal space (resp. 
geometric Poincar\'e complex with Spivak normal structure). According
to the theory of Quinn [14] there is a geometric theory of
transversality for normal spaces, so that by analogy with the
Pontrjagin-Thom isomorphism for smooth bordism there is defined an
isomorphism 
$$ \eqalign{
&\Omega^N_n(X,\nu) \to  \pi_{n+k}(T(\nu))~ ;\cr
&(f:M^n\to X,b:\nu_M\to \nu) \mapsto 
(T(b)(\rho_M):S^{n+k}\,\raise4pt\hbox{$\rho_M \atop \to$}\,T(\nu_M)\,
\raise4pt\hbox{$T(b) \atop \to$}\,T(\nu))~ ,}$$
with inverse 
$$\eqalign{&\pi_{n+k}(T(\nu)) \to \Omega^N_n(X,\nu)~;\cr
&(\rho:S^{n+k}\to T(\nu))\mapsto
(f=\rho\vert:M^n=\rho^{-1}(X)\to X,b:\nu_M\to \nu)~.}$$
The geometric Poincar\'e and normal bordism groups for $n \ge 5$ 
are related by the Levitt-Jones-Quinn exact sequence 
$$\dots \to L_n(\Z[\pi_1(X)]) \to \Omega^P_n(X,\nu) \to \Omega^N_n(X,\nu) \to
L_{n-1}(\Z[\pi_1(X)]) \to \dots~.$$
If $\nu:X \to BG(k)$ admits a $TOP$ reduction $\widetilde{\nu}:X \to BTOP(k)$
the forgetful maps from manifold to normal space bordism
$\Omega_n(X,\nu) \to \Omega^N_n(X,\nu)$
are isomorphisms, and
$$\Omega^P_n(X,\nu) ~=~L_n(\Z[\pi_1(X)])\oplus \Omega^N_n(X,\nu)~.$$

A map of $n$-dimensional normal spaces 
$$(f,b,c) : (M^n,\nu_M,\rho_M)\to (X^n,\nu_X,\rho_X)$$ 
is defined by a map of fibrations $(f,b):(M,\nu_M)\to (X,\nu_X)$ 
together with a homotopy 
$$c~ :~ T(b)\rho_M~ \simeq~ \rho_X~ :~ S^{n+k} \to T(\nu_X)~.$$
The mapping cylinder of $f$ 
$$M(f)~=~M\times [0,1]\cup X/\{(x,1)=f(x)\,\vert\, x\in M\}$$
defines a cobordism $(M(f);M,X)$ of normal spaces, identifying 
$$M~=~M\times\{0\} \subset M(f)~.$$ 
If $M^n$ and $X^n$ are
Poincar\'e complexes the corresponding element of the relative bordism
group is just the surgery obstruction 
$$(M(f);M\cup -X)~=~\sigma_*(f,b)
\in \Omega^{N,P}_{n+1}(X,\nu_X)~=~L_n(\Z[\pi_1(X)])~.$$
Ignoring questions of finite-dimensionality 
(or assuming that $X$ is a finite $n$-dimen\-sional $CW$ complex) 
it is therefore possible to define the inverse isomorphism to
$\Omega^N_n(X,\nu)\to\pi_{n+k}(T(\nu))$ by
$$\pi_{n+k}(T(\nu)) \to \Omega^N_n(X,\nu)~ ;~ \rho \mapsto (X,\nu,\rho)~,$$
without an appeal to the transversality of normal spaces.  The group
$\pi_{n+k}(T(\nu))$ consists of the equivalence classes of normal
structures $(\nu_X:X\to BG(k),\rho_X:S^{n+k}\to T(\nu_X))$ on $X$ with
$\nu_X=\nu$.

\medskip

Following Weiss [19] we shall now identify the algebraic normal bordism
groups $\widehat{L}\langle B,\beta \rangle^n(A)$ with the twisted
quadratic groups $Q_n(B,\beta)$, the algebraic analogues of the
homotopy groups of the Thom space $\pi_{n+k}(T(\nu))$.

\medskip

A {\it cobordism} of $n$-dimensional normal complexes $(C,\phi,\gamma,\theta)$,
$(C',\phi',\gamma',\theta')$ is defined by an $(n+1)$-dimensional symmetric
pair 
$$((f~f'):C\oplus C'\to D,(\delta\phi,\phi\oplus -\phi'))$$
together with bundle maps 
$$(f,\zeta)~ :~ (C,\gamma) \to (D,\delta\gamma)~~,~~
(f',\zeta')~ :~(C',\gamma') \to (D,\delta\gamma)$$
and an equivalence of hyperquadratic structures on $D$ 
$$\delta\theta~ :~ J(\delta\phi) -
(\delta\phi_0;f\phi_0f^*,f'\phi'_0f^{\prime *})^{\%}(S^n\delta\gamma)
+ f^^{
(\phi'_0)^{\%}(S^n\zeta')) \to 0~.$$ 
Similarly for the cobordism of $(B,\beta)$-normal complexes. 
\medskip

The {\it symmetric $(B,\beta)$-structure $L$-groups of $A$}
$L\langle B,\beta\rangle ^n(A)$ $(n\ge 0)$ of
Weiss [19] are the cobordism groups of $n$-dimensional $(B,\beta)$-normal
Poincar\'e complexes over $A$ $(C,\phi,\gamma,\theta,f,\chi)$. 
For the $\cases \text{universal}\cr \text{zero}\endcases$ chain bundle
$\cases (B(\infty),\beta(\infty))\cr (0,0)\endcases$ over $A$ these are just the
$\cases \text{symmetric}\cr \text{quadratic}\endcases$ $L$-groups 
$$\cases L\langle B(\infty),\beta(\infty)\rangle^n(A)~=~L^n(A)\cr
L\langle 0,0\rangle^n(A)~=~L_n(A)~.\endcases$$

The {\it symmetric $(B,\beta)$-structure
$\widehat{L}$-groups}
$\widehat{L}\langle B,\beta \rangle^n(A)$ $(n\ge 0)$
are the cobordism groups of $n$-dimensional $(B,\beta)$-normal complexes
over $A$. 
For the $\cases \text{universal}\cr \text{zero}\endcases$ chain bundle
$\cases (B(\infty),\beta(\infty))\cr (0,0)\endcases$ over $A$ these are just the
$$\cases \widehat{L}\langle B(\infty),\beta(\infty)\rangle^n(A)~=~
\widehat{L}^n(A)\cr
\widehat{L}\langle 0,0\rangle^n(A)~=~0~.\endcases$$

Algebraic surgery was used in Ranicki [15] to prove that
every $n$-dimensional quadratic Poincar\'e complex $(C,\psi)$ is cobordant
to a highly-connected complex $(C',\psi')$, with 
$$H_r(C')~=~0~~(2r\le n-2)~.$$
The boundary of an $n$-dimensional normal
complex $(C,\phi,\gamma,\theta)$ is an $(n-1)$-dimen\-sional quadratic
Poincar\'e complex $(\partial C,\psi)$. Glueing on to
$(C,\phi,\gamma,\theta)$ the trace of the surgery making $(\partial C,\psi)$
highly-connected there is obtained an $n$-dimensional normal complex
$(C',\phi',\gamma',\theta')$ which is cobordant to $(C,\phi,\gamma,\theta)$
and which has a highly-connected boundary, with
$$H_r(\partial C')~=~H_{r+1}(\phi'_0:C^{\prime n-*}\to C')~=~0~(2r\le n-3)~.$$ 
In particular, this shows that every normal
complex is cobordant to a connected complex. Thus
$\widehat{L}\langle B,\beta \rangle^n(A)$ is also the cobordism group of 
connected
$n$-dimensional $(B,\beta)$-normal complexes over $A$. The one-one
correspondence established in \S7 between connected $n$-dimensional
normal complexes and $n$-dimensional (symmetric, quadratic) Poincar\'e
pairs generalizes to a one-one correspondence between connected
$n$-dimensional $(B,\beta)$-normal complexes over $A$ and $n$-dimensional
(symmetric, quadratic) $(B,\beta)$-normal Poincar\'e pairs over $A$, for any
chain bundle $(B,\beta)$ over $A$. It follows that
$\widehat{L}\langle B,\beta \rangle^n(A)$ 
can be identified with the cobordism group
of $n$-dimensional (symmetric, quadratic) $(B,\beta)$-normal Poincar\'e
pairs, and that there is defined an exact sequence 
$$\dots \to L_n(A) \to 
L\langle B,\beta \rangle^n(A)
 \to \widehat{L}\langle B,\beta \rangle^n(A)~
\raise4pt\hbox{$\partial \atop \to$}~L_{n-1}(A) \to \dots~ ,$$
with $\partial$ defined by the quadratic boundary 
$$\partial~:~\widehat{L}\langle B,\beta \rangle^n(A) \to 
L_{n-1}(A)~ ;~ (C,\phi,\gamma,\theta,f,\chi) \mapsto 
\partial (C,\phi,\gamma,\theta)~.$$ 

A map of $n$-dimensional normal complexes 
$$(f,\xi,\chi,\eta)~ :~ (C,\phi,\gamma,\theta) \to 
(C',\phi',\gamma',\theta')$$
determines an $(n+1)$-dimensional symmetric
pair $((f~1):C\oplus C'\to C',(\xi,\phi\oplus -\phi'))$, bundle maps
$$(f,\chi)~ :~ (C,\gamma) \to (C',\gamma')~,~
(1,0)~ :~ (C',\gamma') \to (C',\gamma')$$
and an equivalence of hyperquadratic structures on $C'$
$$\eta : J(\xi) - (\xi_0;f\phi_0f^*,\phi'_0)^{\%}(S^n\gamma') 
+ f^{\%}(\theta - (\phi_0)^{\% }(S^n\gamma)) -\theta' \to 0~ ,$$ 
defining a cobordism between $(C,\phi,\gamma,\theta)$
and $(C',\phi',\gamma',\theta')$ by analogy with the mapping cylinder
construction of geometric normal bordisms. Similarly for maps of
$(B,\beta)$-normal complexes. 
It follows that the abelian group morphisms 
$$\eqalign{&\widehat{L}\langle B,\beta \rangle^n(A) \to Q_n(B,\beta)~;~
(C,\phi,\gamma,\theta,f,\chi) \mapsto (f,\chi)_{\% }(\phi,\theta)~,\cr
&Q_n(B,\beta) \to \widehat{L}\langle B,\beta \rangle^n(A)~;~
(\phi,\theta) \mapsto (B,\phi,\beta,\theta,1,0)}$$ 
are inverse isomorphisms.
\medskip

For example, if $A=\Z$ and $(B(\infty),\beta(\infty))$ is the universal
chain bundle over $\Z$ (as constructed at the end of \S6) then
$$\eqalign{
&L\langle B(\infty),\beta(\infty) \rangle^n(\Z)~=~L^n(\Z)~=~
\cases \Z &\text{if $n \equiv 0(\bmod\,4)$}\cr
\Z_2&\text{if $n\equiv 1(\bmod\,4)$}\cr
0&\text{if $n\equiv 2,3(\bmod\,4)$~,}\endcases\cr
&L_n(\Z)~=~\cases \Z &\text{if $n \equiv 0(\bmod\,4)$}\cr
\Z_2&\text{if $n\equiv 2(\bmod\,4)$}\cr
0&\text{if $n\equiv 1,3(\bmod\,4)$~,}\endcases\cr
&\widehat{L}\langle B(\infty),\beta(\infty) \rangle^n(\Z)~=~
Q_n(B(\infty),\beta(\infty))~
=~\cases \Z_8&\text{if $n \equiv 0(\bmod\,4)$}\cr
\Z_2&\text{if $n\equiv 1,3(\bmod\,4)$}\cr
0&\text{if $n\equiv 2(\bmod\,4)$~.}\endcases}$$

\medskip

A spherical fibration $\nu:X \to BG(k)$ determines a chain bundle 
$(C(\widetilde{X}),\gamma)$ over $\Z[\pi_1(X)]$ ([15], [19])
and there is defined a natural transformation of exact
sequences from the Levitt-Jones-Quinn Poincar\'e bordism sequence

$$\xymatrix@C-10pt{\dots \ar[r] 
& L_n(\Z[\pi_1(X)]) \ar@{=}[d] \ar[r]
& \Omega^P_n(X,\nu) \ar[r] \ar[d]
& \pi_{n+k}(T(\nu)) \ar[r] \ar[d]
& L_{n-1}(\Z[\pi_1(X)]) \ar@{=}[d] \ar[r]&\dots \\
\dots \ar[r]
& L_n(\Z[\pi_1(X)])  \ar[r] 
& L^n(C(\widetilde X),\gamma) \ar[r] 
& Q_n(C(\widetilde X),\gamma) \ar[r]
& L_{n-1}(\Z[\pi_1(X)])  \ar[r]&\dots}$$

\medskip

\noindent with $\Omega^P_n(X,\nu)  \to L^n(C(\widetilde X),\gamma)$
a generalized symmetric signature map.

\medskip

\noindent \S9. {\bf Normal Wu classes}
\medskip

The Wu classes of the symmetric structure $\phi$ and the bundles
$\beta,\gamma$ in an $n$-dimensional $(B,\beta)$-normal complex
$(C,\phi,\gamma,\theta,f,\chi)$ are related by a commutative diagram

$$\xymatrix{H^{n-r}(C) \ar[r]^-{\di{\phi_0}} \ar[d]_-{\di{v_r(\phi)}} & 
H_r(C) \ar[r]^-{\di{f_*}} \ar[d]_-{\di{\hbox{$\widehat{v}_r(\gamma)$}}}
& H_r(B) \ar[dl]^-{\di{\hbox{$\widehat{v}_r(\beta)$}}}\\
H^{n-2r}(\Z_2;A,(-1)^{n-r}) \ar[r]^-{\di{J}} & \hbox{$\widehat{H}^r(\Z_2;A)$} 
&}$$

\medskip

\noindent For any chain bundle $(B,\beta)$ and any chain complex $C$ 
we shall now define symmetric $(B,\beta)$-structure groups 
$Q\langle B,\beta \rangle^n(C)$ ($n \ge 0$) to fit into an exact sequence
$$\dots \to Q\langle B,\beta \rangle^n(C)\to Q^n(C)\oplus H_n(B^t\otimes_AC)
 \to \widehat{Q}^n(C)\to Q\langle B,\beta \rangle ^{n-1}(C) \to \dots~.$$
The Wu classes $v_r(\phi)$ of a symmetric complex $(C,\phi)$ will
then be refined to the normal Wu classes of a $(B,\beta)$-normal complex
$(C,\phi,\gamma,\theta,f,\chi)$ 
$$v_r~=~v_r(\phi,\gamma,\theta,f,\chi)~ :~
H^{n-r}(C) \to Q\langle B,\beta \rangle^n(S^{n-r}A)~,$$
with
$$v_r(\phi)~:~H^{n-r}(C) \xymatrix@C-10pt{\ar[r]^-{\di{v_r}}&}
Q\langle B,\beta \rangle^n(S^{n-r}A) 
\to Q^n(S^{n-r}A)~=~H^{n-2r}(\Z_2;A,(-1)^{n-r})~.$$
In \S11 below the normal Wu
classes will be used to define a $\Z_4$-valued quadratic function on
$H^n(C)$ for a $2n$-dimensional symmetric Poincar\'e complex $(C,\phi)$ over
$\Z_2$ with normal $(v_{n+1}=0)$-structure, as required to define the
$\Z_8$-valued invariant of Brown [4]. 
\medskip

Let $(B,\beta)$ be a chain bundle over $A$, and let $C$ be a finite-dimensional $A$-module chain
complex. An {\it $n$-dimensional symmetric $(B,\beta)$-structure on $C$}
$(\phi,\theta,f)$ is defined by an $n$-dimensional symmetric structure
$\phi\in (W^{\% }C)_n$ together with a chain 
$\theta\in (\widehat{W}^{\%}C)_{n+1}$ and a chain map $f:B^{n-*}\to C$ such
that 
$$J(\phi) - \widehat{f}^{\% }(S^n\beta)~=~d\theta \in
(\widehat{W}^{\% }C)_n~.$$
An $n$-dimensional $(B,\beta)$-normal structure
$(\phi,\gamma,\theta,g,\chi)$ on $C$ determines the $n$-dimen\-sional symmetric
$(B,\beta)$-structure 
$(\phi,\theta+(\phi_0)^{\% }(S^n\chi),\phi_0g^*)$
on $C$. Conversely, if $f^*:C^{n-*}\to B$ is a composite 
$$f^* ~:~C^{n-*}~\xymatrix@C20pt{\ar[r]^-{\di{\phi_0}}&}~C~
\xymatrix@C20pt{\ar[r]^-{\di{g}}&}~B$$ 
(as is always the case up to chain homotopy if $(C,\phi)$ is a 
Poincar\'e complex) the symmetric $(B,\beta)$-structure $(\phi,\theta,f)$ 
determines the $n$-dimensional $(B,\beta)$-normal structure 
$(\phi,g^*\gamma,\theta,g,0)$. 
\medskip

An {\it $n$-dimensional symmetric $(B,\beta)$-structure
{\rm(}Poincar\'e{\rm)} complex over $A$} $(C,\phi,\allowbreak\theta,f)$ is an
$n$-dimensional $A$-module chain complex $C$ together with an
$n$-dimensional symmetric $(B,\beta)$-structure $(\phi,\theta,f)$ (such
that $\phi_0:C^{n-*}\to C$ is a chain equivalence).  As for symmetric
(Poincar\'e) pairs there is also the analogous notion of {\it symmetric
$(B,\beta)$-structure {\rm(}Poincar\'e{\rm)} pair}.  There is
essentially no difference between symmetric $(B,\beta)$-structure
Poincar\'e complexes and $(B,\beta)$-normal Poincar\'e complexes, so
that the $L$-groups $L\langle B,\beta \rangle^n(A)$ $(n\ge 0)$ can also
be regarded as the cobordism groups of $n$-dimensional symmetric
$(B,\beta)$-structure Poincar\'e complexes over $A$.
\medskip

An {\it equivalence} of $n$-dimensional symmetric
$(B,\beta)$-structures on $C$ 
$$(\xi,\eta,g)~ :~ (\phi,\theta,f) \to (\phi',\theta',f')$$ 
is defined by an equivalence of symmetric
structures $\xi:\phi\to \phi'$ together with a chain $\eta\in
(\widehat{W}^{\% }C)_{n+2}$ and a chain homotopy $g:f\simeq f':B^{n-*}\to
C$ such that 
$$J(\xi) - (g;f,f')^{\% }(S^n\beta) - \theta' +
\theta~=~d\eta \in (\widehat{W}^{\% }C)_{n+1}~.$$ 

The {\it $n$-dimensional symmetric $(B,\beta)$-structure group of $C$}
$Q\langle B,\beta \rangle^n(C)$
is the abelian group of equivalence classes of $n$-dimensional symmetric
$(B,\beta)$-symmetric structures on $C$, with addition by 
$$(\phi,\theta,f) +
(\phi',\theta',f')~=~(\phi+\phi',\theta+\theta'+[f,f'](S^n\beta),f+f')
\in Q\langle B,\beta \rangle^n(C)~.$$
There is also a more economical description of
$Q\langle B,\beta \rangle^n(C)$ 
as the abelian group of equivalence classes of pairs
$(\psi,f)$ defined by an $n$-dimensional quadratic structure 
$\psi\in (W_{\%}C)_n$ and a chain map $f:B^{n-*}\to C$ such that 
$$f_{\% }H(S^n\beta)~=~d\psi \in (W_{\% }C)_{n-1}~,$$
so that up to signs
$$f\beta_{-n-s-1}f^*~=~d\psi_s + \psi_sd^* + \psi_{s+1} +
\psi^*_{s+1} \in \Hom_A(C^{-*},C)_{n-s}~~(s\ge 0)~,$$ 
subject to the equivalence relation
\medskip

{\parskip=2pt
\parindent=100pt
\item{$(\psi,f) \sim (\psi',f')$} if there exist a chain homotopy 
$g:f\simeq f':B^{n-*}\to C$ 
\item{} and an equivalence of quadratic structures 
\item{} $\chi~:~\psi'-\psi \to (g;f,f')_{\% }H(S^n\beta)$~,}

\noindent with addition by 
$$(\psi,f) + (\psi',f')~=~(\psi+\psi'+H([f,f'](S^n\beta)),f+f')~.$$ 
\medskip

\noindent The pair $(\psi,f)$ determines the triple $(\phi,\theta,f)$ with 
$$\eqalign{
&\phi_s~=~\cases f\beta_{s-n}f^*&\text{if $s\ge 1$}\cr
f\beta_{-n}f^*+(1+T)\psi_0&\text{if $s=0$~,}\endcases\cr
&\theta_s~=~\cases 0&\text{if $s\ge 0$}\cr \psi_{-s-1}&\text{if $s\le -1$~.}
\endcases}$$
Conversely, a triple $(\phi,\theta,f)$ determines the pair $(\psi,f)$ with 
$$\psi_s~=~\theta_{-s-1}~~(s\ge 0)~.$$ 

Given an $n$-dimensional symmetric
$(B,\beta)$-structure $(\phi,\theta,f)$ on $C$, a chain bundle map 
$(g,\chi) : (B,\beta)\to (B',\beta')$ and a chain map 
$h: C\to C'$ define the
{\it pushforward} $n$-dimensional symmetric $(B',\beta')$-structure on $C'$
$$\langle g,\chi\rangle(h)^{\% }(\phi,\theta,f)~=~
(h^{\% }(\phi),h^{\%}(\theta+S^n(\widehat{f}^{\% }\chi)),hfg^*)~.$$ 
Thus the groups
$Q\langle B,\beta \rangle^*(C)$ are covariant in both $(B,\beta)$ and $C$, 
with pushforward abelian group morphisms 
$$\eqalign{\langle g,\chi \rangle (h)^{\%}~ :~~&
Q\langle B,\beta \rangle^n(C) \to Q\langle B',\beta'\rangle^n(C')~;\cr
&(\phi,\theta,f) \mapsto (h^{\%}(\phi),\widehat{h}^{\%}(
\theta+\widehat{f}^{\%}(S^n\chi)),hfg^*)}$$
depending only on the homotopy classes of $(g,\chi)$ and $h$. 
\medskip

An $n$-dimensional symmetric $(B,\beta)$-structure $(\phi,\theta,f)$ on $C$
determines an $n$-dimen\-sional symmetric structure $\phi\in (W^{\% }C)_n$ on
$C$, so that there is defined a forgetful map 
$$s~:~Q\langle B,\beta \rangle^n(C) \to Q^n(C)~;~
 (\phi,\theta,f) \mapsto \phi~ .$$ 
 An $n$-dimensional quadratic structure $\psi\in (W_{\% }C)_n$ on $C$ 
 determines an $n$-dimension\-al symmetric $(B,\beta)$-structure 
 $((1+T)\psi,\theta,0)$ on $C$ for any
$(B,\beta)$, with 
$$\theta_s~=~\cases \psi_{-s-1}&\text{if $s\le -1$}\cr
0&\text{if $s\ge 0$~.}\endcases$$
Thus there are also defined forgetful maps 
$$s~ :~ Q_n(C) \to Q\langle B,\beta \rangle^n(C)~ ;~ \psi \mapsto 
((1+T)\psi,\theta,0)~,$$
and $sr=1+T:Q_n(C)\to Q^n(C)$.
\medskip

Let $P\langle B,\beta \rangle^n(C)$ 
be the abelian group of equivalence classes of
$n$-dimensional symmetric $(B,\beta)$-structures $(\phi,\theta,f)$ with
$\phi=0$, to be denoted $(\theta,f)$, subject to the equivalence relation
\medskip

{\parskip=2pt
\parindent=100pt
\item{$(\theta,f) \sim (\theta',f')$} if there exists an equivalence of
$(B,\beta)$-structures 
\item{} $(0,\eta,g) : (0,\theta,f)\to (0,\theta',f')$.\par}
\medskip

\noindent 
The symmetric $(B,\beta)$-structure groups $Q\langle B,\beta \rangle^*(C)$ 
and the groups $P\langle B,\beta \rangle^*(C)$ are related by a 
commutative braid of exact sequences of abelian groups 

$$\xymatrix@C-20pt{
Q_n(C)
\ar[dr]\ar@/^2pc/[rr]^{\di{1+T}}&&
Q^n(C)\ar[dr]^{\di{J}} \ar@/^2pc/[rr]^{} &&P\langle B,\beta \rangle^{n-1}(C) \\&
Q\langle B,\beta \rangle^n(C)\ar[ur] \ar[dr] && 
\hbox{$\widehat{Q}^n(C)$} \ar[ur] \ar[dr]^{\di{H}}&&\\
P\langle B,\beta \rangle^n(C)
\ar[ur]\ar@/_2pc/[rr]_-{}&&H_n(B^t\otimes_AC)
\ar[ur]\ar@/_2pc/[rr]_{}&&Q_{n-1}(C)
}$$

\medskip

\noindent If $(B,\beta)$ is the 
$\cases \text{universal}\cr \text{zero}\endcases$ chain bundle
$\cases (B(\infty),\beta(\infty))\cr (0,0)\endcases$ the forgetful map
$$\cases 
Q\langle B(\infty),\beta(\infty)\rangle^n(C) \to Q^n(C)~;~
(\phi,\theta,f) \mapsto \phi\cr
Q_n(C) \to Q\langle 0,0\rangle^n(C)~;~\psi \mapsto ((1+T)\psi,\theta)\endcases$$
is an isomorphism and 
$$\cases P\langle B(\infty),\beta(\infty)\rangle^n(C)~=~0\cr
P\langle 0,0 \rangle^n (C)~=~\widehat{Q}^{n+1}(C)~.\endcases$$

The {\it Wu classes} of an $n$-dimensional symmetric
$(B,\beta)$-structure $(\phi,\theta,f)$ on $C$ are the $A$-morphisms 
$$\eqalign{
v_r(\phi,\theta,f)~:~&H^{n-r}(C) \to Q\langle B,\beta \rangle^n(S^{n-r}A)~;\cr
&(x:C\to S^{n-r}A) \mapsto \langle 1,0\rangle (x)^{\%}(\phi,\theta,f)~.}$$
Now 
$$\eqalign{&Q\langle B,\beta \rangle^n(S^{n-r}A)~=\cr
&\cases H_r(B)\! &\text{if $2r<n$}\cr
\{(a,b)\in A\oplus B_r\vert 
db=0\in B_{r-1}\}/\!\!\sim \!&\text{if $2r=n$}\cr
\{(a,b)\in A\oplus B_r\vert 
a+(-1)^r\overline{a}+\beta_{-2r}(b)(b)=0\in A,db=0\in
B_{r-1}\}/\sim \!&\text{if $2r>n$}\endcases}$$
with the equivalence relation~$\sim$~defined by 
\medskip

{\parskip=2pt
\parindent=100pt
\item{$(a,b) \sim (a',b')$} if there exists $(x,y)\in A\oplus B_{r+1}$ such that 
\item{} $a' - a = x +(-1)^{r+1}\overline{x} + \beta_{-2r-2}(y)(y)$
\item{} \hskip75pt $+\beta_{-2r-1}(y)(b) +\beta_{-2r-1}(b')(y) \in A$,
\item{} $b'-b=dy \in B_r$\par}
\medskip

\noindent and addition by 
$$(a,b) + (a',b')~=~(a+a'+\beta_{-2r}(b)(b'),b+b')~.$$ 
The map to the symmetric $Q$-group is given by
$$\eqalign{&
Q\langle B,\beta \rangle^n(S^{n-r}A) \to Q^n(S^{n-r}A)~=~
H^{n-2r}(\Z_2;A,(-1)^{n-r})~;\cr
&\hskip75pt 
\cases b \mapsto \beta_{-2r}(b)(b)&\text{if $2r<n$}\cr
(a,b) \mapsto a+(-1)^r\overline{a}+\beta_{-2r}(b)(b)&\text{if $2r=n$}\cr
0&\text{if $2r>n$~.}\endcases}$$
The Wu classes are given by 
$$\eqalign{v_r(\phi,\theta,f)~:~&H^{n-r}(C) \to 
Q\langle B,\beta \rangle^n(S^{n-r}A)~;\cr
&z \mapsto 
\cases f^*(z) &\text{if $2r<n$} \cr
(\theta_{n-2r-1}(z)(z),f^*(z))&\text{if $2r\ge n$}\endcases~~
(z\in C^{n-r},d^*z=0) ~.}$$

\noindent \S10. {\bf Forms}
\medskip

In Ranicki [15] the even-dimensional
$\cases \text{symmetric}\cr \text{quadratic} \endcases$ $L$-groups
$\cases L^{2n}(A)\cr L_{2n}(A)\endcases$ $(n\ge 0)$ were related
to the Witt groups $\cases W^{(-1)^n}(A)\cr W_{(-1)^n}(A)\endcases$
of nonsingular $\cases \text{$(-1)^n$-symmetric}\cr
\text{$(-1)^n$-quadratic}\endcases$ forms over $A$.
In particular, it was shown that
$$\cases L^0(A)~=~W^{+1}(A)\cr L_{2n}(A)~=~W_{(-1)^n}(A)\endcases~~(n\ge 0)~.$$
This relationship between $L$-groups and Witt groups will now be
generalized to the even-dimensional symmetric $(B,\beta)$-structure
$L$-groups $L\langle B,\beta \rangle^{2n}(A)$ and the Witt groups 
$W_{Q(n)}(A)$ of nonsingular $Q(n)$-quadratic forms over $A$, with 
$(B,\beta)$ any chain bundle over $A$ and
$$Q(n)~=~Q\langle B,\beta \rangle^{2n}(S^nA)~.$$

Let $\epsilon = \pm 1$.
An {\it $\epsilon$-symmetric form over $A$}
$(M,\lambda)$ is a f.g. projective $A$-module $M$ together with an element
$\lambda\in \Hom_A(M,M^*)$ such that 
$$\epsilon \lambda^*~=~\lambda~ :~ M \to M^*~ .$$ 
Equivalently, the form is defined by a pairing 
$$\lambda~ :~ M \times M \to A~ ;~ (x,y) \to 
\lambda(x,y)~=~\lambda(x)(y)$$ such that
$$\eqalign{&\lambda(ax,by)~=~b\lambda(x,y)\overline{a}~ ,\cr
&\lambda(x+x',y)~=~\lambda(x,y) + \lambda(x',y)~,\cr
&\epsilon \overline{\lambda(y,x)}~=~\lambda(x,y)~~
(x,y\in M, a,b\in A)~.}$$
Let $Q(\epsilon )$ be an $A$-group together with $A$-morphisms 
$$\eqalign{
&r~:~Q(\epsilon ) \to 
H^0(\Z_2;A,\epsilon )~=~\{a\in A\,\vert\, \epsilon \overline{a}=a\}~,\cr
&s~:~
H_0(\Z_2;A,\epsilon )~=~
A/\{b-\epsilon \overline{b}\,\vert\, b\in A\} \to Q(\epsilon )}$$
such that 
$$rs~=~1+T_{\epsilon }~:~H_0(\Z_2;A,\epsilon ) \to H^0(\Z_2;A,\epsilon )~.$$ 
A {\it $Q(\epsilon )$-quadratic form over $A$}
$(M,\lambda,\mu)$ is an $\epsilon$-symmetric
form $(M,\lambda)$ together with an $A$-morphism 
$\mu : M\to Q(\epsilon )$ such that 
$$\eqalign{
&r(\mu(x))~=~\lambda(x,x) \in H^0(\Z_2;A,\epsilon )~,\cr
&\mu(x+y) - \mu(x)
- \mu(y)~=~s(\lambda(x,y)) \in Q(\epsilon )~~(x,y\in M)~.}$$
There is an
evident notion of {\it isomorphism} of $Q(\epsilon )$-quadratic forms. 
\medskip

A $Q(\epsilon )$-quadratic form $(M,\lambda,\mu)$ is {\it nonsingular} if
$\lambda\in \Hom_A(M,M^*)$ is an isomorphism of $A$-modules. 
\medskip

A nonsingular $Q(\epsilon )$-quadratic form $(M,\lambda,\mu)$ is 
{\it hyperbolic} if there exists a direct summand $L\subset M$ such that
\medskip

{\parskip=2pt
\parindent=23pt
\item{(i)} the inclusion $j\in \Hom_A(L,M)$ fits into an exact sequence 
$$0 \to L~\raise4pt\hbox{$j\atop\to$}~M~
\raise4pt\hbox{$j^*\lambda\atop\to$}~L^* \to 0~,$$ 
\item{(ii)} $\mu j~=~0~ :~ L \to Q(\epsilon )~.$\par}
\medskip

The {\it $Q(\epsilon)$-quadratic Witt group of $A$} 
$W_{Q(\epsilon )}(A)$ is the abelian group of
equivalence classes of nonsingular $Q(\epsilon )$-quadratic forms
$(M,\lambda,\mu)$, subject to the equivalence relation
\medskip

{\parskip=2pt
\parindent=125pt
\item{$(M,\lambda,\mu) \sim (M',\lambda',\mu')$} if there exists an
isomorphism
\item{} 
$(M,\lambda,\mu)\oplus (N,\nu,\rho) \to (M',\lambda',\mu')\oplus (N',\nu',\rho')$
\item{}
for some hyperbolic $Q(\epsilon )$-quadratic forms
\item{} $(N,\nu,\rho)$, $(N',\nu',\rho')$.\par}
\medskip

For $Q(\epsilon )=H^0(\Z_2;A,\epsilon )$, $r=1$, $s=1+T_{\epsilon}$ a 
$Q(\epsilon )$-quadratic form $(M,\lambda,\mu)$ may be
identified with the $\epsilon$-symmetric form $(M,\lambda)$, since $\lambda$
determines $\mu$ by 
$$\mu(x)~=~\lambda(x,x) \in H^0(\Z_2;A,\epsilon )~~(x \in M)~.$$ 
The Witt group of $\epsilon$-symmetric forms $W_{Q(\epsilon )}(A)$ 
is denoted by $W^{\epsilon}(A)$. 
\medskip

For $Q(\epsilon )=H_0(\Z_2;A,\epsilon )$, $r=1+T_{\epsilon }$, $s=1$ a
$Q(\epsilon )$-quadratic form $(M,\lambda,\mu)$ is just a $\epsilon$-quadratic
form in the sense of Wall [18]. The Witt group of $\epsilon$-quadratic forms 
$W_{Q(\epsilon)}(A)$ is denoted by $W_{\epsilon }(A)$. 
\medskip

For $Q(\epsilon )=\im
(1+T_{\epsilon}:H_0(\Z_2;A,\epsilon )\to H^0(\Z_2;A,\epsilon ))$, 
$r=\text{projection}$, $s=\text{injection}$ a
$Q(\epsilon )$-quadratic form $(M,\lambda,\mu)$ is just an $\epsilon$-symmetric 
form $(M,\lambda)$ for which there exists an $\epsilon$-quadratic form
$(M,\lambda,\mu:M\to H_0(\Z_2;A,\epsilon ))$. 
Such an $\epsilon$-symmetric form is {\it even}. 
The Witt group of even $\epsilon$-symmetric forms $W_{Q(\epsilon )}(A)$
is denoted by $W\langle v_0 \rangle^{\epsilon }(A)$. 
\medskip

For $\epsilon =+1$
$\cases \text{$\epsilon$-symmetric}\cr
\text{$\epsilon$-quadratic}\endcases$ is abbreviated to
$\cases \text{symmetric}\cr \text{quadratic}~.\endcases$ 
\medskip

A $2n$-dimensional $\cases \text{symmetric}\cr \text{quadratic}\endcases$
(Poincar\'e) complex over $A$
$\cases (C,\phi)\cr (C,\psi)\endcases$ with $H^n(C)$ a f.g. projective
$A$-module determines a (nonsingular)
$\cases \text{$(-1)^n$-symmetric}\cr
\text{$(-1)^n$-quadratic}\endcases$ form over $A$
$\cases (H^n(C),\phi_0,v_n(\phi))\cr (H^n(C),(1+T)\psi_0,v^n(\psi))\endcases$
with 
$$\cases 
v_n(\phi)~:~H^n(C) \to H^0(\Z_2;A,(-1)^n)~;~x \mapsto \phi_0(x)(x)\cr
v^n(\psi)~ :~ H^n(C) \to H_0(\Z_2;A,(-1)^n)~;~x \mapsto \psi_0(x)(x)~.\endcases$$ 
Conversely, a (nonsingular) 
$\cases \text{$(-1)^n$-symmetric}\cr \text{$(-1)^n$-quadratic}\endcases$ form
$\cases (M,\lambda)\cr (M,\lambda,\mu)\endcases$ determines a
$2n$-dimensional $\cases \text{symmetric}\cr \text{quadratic}\endcases$ 
(Poincar\'e) complex $\cases (C,\phi)\cr (C,\psi)\endcases$ such that 
$$\eqalign{&
\cases \phi_0\cr (1+T)\psi_0\endcases~=~\lambda~:~C^n~=~M \to C_n~=~M^*~~,~~
C_r~=~0~~ (r\ne n)~,\cr
&v^n(\psi)~=~\mu~:~H^n(C)~=~M \to H_0(\Z_2;A,(-1)^n)~.}$$
The corresponding morphisms from the Witt groups to the $L$-groups 
$$\cases W^{(-1)^n}(A) \to L^{2n}(A)~;~(M,\lambda) \mapsto (C,\phi)\cr
W_{(-1)^n}(A) \to L_{2n}(A)~;~(M,\lambda,\mu) \mapsto (C,\psi)\endcases$$
were shown in Ranicki [15] to be isomorphisms for $n=0$ 
if $A$ is any ring, and for all $n \ge 0$ 
if $A$ is $\cases \text{a Dedekind}\cr \text{any}\endcases$ ring. 
For a Dedekind ring $A$ the inverse isomorphism in symmetric $L$-theory is
given by 
$$L^{2n}(A) \to W^{(-1)^n}(A)~ ;~(C,\phi) \mapsto 
(H^n(C)/(\text{torsion}),\phi_0)~.$$ 
The inverse isomorphism in quadratic $L$-theory is given for any $A$ by
$$\eqalign{&L_{2n}(A) \to W_{(-1)^n}(A)~;\cr
&(C,\psi) \mapsto (\coker(\pmatrix d^* & 0 \cr
(1+T)\psi_0 & d \endpmatrix :
C^{n-1} \oplus C_{n+2} \to C^n \oplus C_{n+1}), \bigg[
\matrix \psi_0 & d \cr 0 & 0 \endmatrix \bigg])~.}$$
If $A$ is a field this isomorphism can also be expressed as
$$(C,\psi) \mapsto (H^n(C),(1+T)\psi_0,v^n(\psi))$$ 
but this is not the case in general -- see Milgram and Ranicki [12,\t p.406].

\medskip

Given $A$-groups $M,N$ and a symmetric bilinear pairing 
$$\phi~ :~ N\times N \to M$$ 
such that 
$$\phi(ay,ay')~=~a\phi(y,y') \in M~~ (a\in A, y,y'\in N)$$ 
let $M\times _{\phi}N$ be the $A$-group of pairs $(x\in M,y\in N)$,
with addition by 
$$(x,y) + (x',y')~=~(x+x'+\phi(y,y'),y+y') \in M\times _{\phi}N$$
and $A$ acting by 
$$A \times (M\times _{\phi}N) \to M\times _{\phi}N~ ;~ (a,(x,y)) \mapsto (ax,ay)~ .$$
There is then defined a short exact
sequence of $A$-groups and $A$-morphisms 
$$0 \to M \to M\times _{\phi}N \to N \to 0$$
with
$$ \eqalign{&M \to M\times _{\phi}N~;~x \mapsto (x,0)\cr
&M\times _{\phi}N \to N~ ;~(x,y) \mapsto y~ .}$$

Given a chain bundle $(B,\beta)$ over $A$ define the $A$-group 
$$\eqalign{Q(n)~&=~Q\langle B,\beta \rangle^{2n}(S^nA)\cr
&=~\{(a,b) \in A\oplus B_n\,\vert\, db=0 \in B_{n-1}\}/\!\!\sim~,}$$
where
\medskip

{\parskip=3pt
\parindent=75pt
\item{$(a,b) \sim (a',b')$}
if there exist $(x,y)\in A\oplus B_{n+1}$ such that 
\item{} $a' - a~=~x + (-1)^{n+1}\overline{x} + \beta_{-2n-2}(y)(y)$
\item{} \hphantom{$a' - a~=~x + (-1)^{n+1}\overline{x}$\,} 
$+\,\beta_{-2n-1}(y)(b) +\beta_{-2n-1}(b')(y)~,$
\item{} $b' - b~=~dy~,$\par}
\smallskip

\noindent with addition by 
$$(a,b) + (a',b')~=~(a+a'+\beta_{-2n}(b)(b'),b+b')$$ 
and $A$-action by 
$$A \times Q(n) \to Q(n)~;~(x,(a,b)) \mapsto 
x(a,b)~=~(xa\overline{x},xb)~.$$ 
The $A$-morphisms 
$$\eqalign{
&r~:~Q_{2n}(S^nA)~=~H_0(\Z_2;A,(-1)^n) \to Q(n)~;~a \mapsto (a,0)~,\cr
&s~:~Q(n) \to Q^{2n}(S^nA)~=~H^0(\Z_2;A,(-1)^n)~;\cr
&\hskip100pt (a,b) \mapsto a+(-1)^n\overline{a}+\beta_{-2n}(b)(b)}$$
are such that there is defined a commutative braid of exact sequences

$$\xymatrix@!C@C-45pt{
H_{n+1}(B) \ar[dr]\ar@/^2pc/[rr]^{} && 
Q_{2n}(S^nA) \ar[dr]^-{\di{r}} \ar@/^2pc/[rr]^-{\di{1+T}}&&
Q^{2n}(S^nA) \ar[dr]^-{\di{J}} \ar@/^2pc/[rr]^{}&&
0 \\
& \hbox{$\widehat{Q}^{2n+1}(S^nA)$} \ar[ur]^-{\di{H}} \ar[dr] &&
Q(n) \ar[ur]^-{\di{s}} \ar[dr] && \hbox{$\widehat{Q}^{2n}(S^nA)$} \ar[ur] \ar[dr] & \\
0 \ar[ur]\ar@/_2pc/[rr]_{} && 
P\langle B,\beta \rangle^{2n}(S^nA) 
\ar[ur]\ar@/_2pc/[rr]_{} &&
H_n(B) \ar[ur]\ar@/_2pc/[rr]_{} && 0 
}$$
\medskip

\noindent with
$$Q(n) \to H_n(B)~;~(a,b) \mapsto b~.$$
If $(B,\beta)$ is such that for all $y\in B_{n+1}$ there exists
$x\in A$ such that 
$$(\beta_{-2n-2}+\beta_{-2n-1}d)(y)(y)~=~x +
(-1)^n\overline{x} \in A$$ 
(e.g. if $d=0:B_{n+1}\to B_n$ and 
$v_{n+1}(\beta)=0:H_{n+1}(B)\to \widehat{H}^{n+1}(\Z_2;A)$\t) 
then there is a natural identification of $A$-groups 
$$Q(n)~=~Q_{2n}(S^nA)\times_{\beta_{-2n}}H_n(B)$$ 
with 
$$\beta_{-2n}~:~H_n(B) \times H_n(B) \to Q_{2n}(S^nA)~;~
(b,b')\mapsto \beta_{-2n}(b)(b')~.$$ 
For any chain bundle $(B,\beta)$ and any $Q(n)$-quadratic form 
$(M,\lambda,\mu)$ there exist $A$-module morphisms 
$$g ~:~ M \to B_n~ ~,~ ~ \psi~ :~ M \to M^*$$ 
such that 
$$\eqalign{&dg~=~0~ :~M \to B_{n-1}~ ,\cr
&\lambda - g^*\beta_{-2n}g~=~\psi+(-1)^n\psi^*~:~M \to M^*~,\cr
&\mu~ :~ M \to Q(n)~ ;~ x \mapsto (\psi(x)(x),g(x))~.}$$
If $(M,\lambda,\mu)$ is a nonsingular form there is thus defined a
$2n$-dimensional symmetric $(B,\beta)$-structure Poincar\'e complex
$(C,\phi,\theta,f)$ with 
$$\eqalign{&\phi_0~=~\lambda~ :~ C^n~=~M\to C_n~=~M^*~~,~~
C_r~=~0~~~ (r\ne n)~,\cr
&\theta_{-1}~=~\psi~ :~ C^n~=~M \to C_n~=~M^*~ ,\cr
&f~=~g\lambda^{-1}~:~C_n~=~M^*~\raise4pt\hbox{$\lambda^{-1}\atop\to$}~
M~\raise4pt\hbox{$g\atop \to$}~B_n~,\cr
&v_n(\phi,\theta,f)~=~\mu~:~H^n(C)~=~M \to Q(n)~.}$$
The construction defines a morphism of abelian groups
$$W_{Q(n)}(A) \to L\langle B,\beta \rangle^{2n}(A)~;~
(M,\lambda,\mu) \mapsto (C,\phi,\theta,f)~.$$
Conversely, if $(C,\phi,\theta,f)$ is a $2n$-dimensional symmetric
$(B,\beta)$-structure Poincar\'e complex such that $H^n(C)$ is a f.g. 
projective $A$-module there is defined a nonsingular $Q(n)$-quadratic form
$(H^n(C),\phi_0,v_n(\phi,\theta,f))$, with 
$$v_n(\phi,\theta,f)~ :~
H^n(C) \to Q(n)~;~x \mapsto (\theta_{-1}(x)(x),f(x))~.$$
It follows that for any ring $A$ there is a natural identification of the 
0-dimensional $L$-group with the Witt group 
$$L\langle B,\beta \rangle^0(A)~=~W_{Q(0)}(A)~.$$
For a field $A$ the morphisms 
$$W_{Q(n)}(A) \to L\langle B,\beta \rangle^{2n}(A)~~(n\ge 0)$$
are injections, which are split by 
$$L\langle B,\beta \rangle^{2n}(A) \to W_{Q(n)}(A)~;~
(C,\phi,\theta,f) \mapsto (H^n(C),\phi_0,v_n(\phi,\theta,f))~~(n\ge 0)~.$$

For any ring with involution $A$ let $(B(\infty),\beta(\infty))$ be
the universal chain bundle of Weiss [19] (cf. \S6 above), with 
isomorphisms 
$$\eqalign{
&\widehat{v}_m(\beta(\infty))~:~H_m(B(\infty)) \to \widehat{H}^m(\Z_2;A)~,\cr
&L\langle B(\infty),\beta(\infty) \rangle^m(A)~\cong~L^m(A)}$$
and an exact sequence
$$\dots \to L_m(A)~\raise4pt\hbox{$1+T\atop \to$}~L^m(A) \to 
Q_m(B(\infty),\beta(\infty))~
\raise4pt\hbox{$\partial\atop\to$}~L_{m-1}(A) \to \dots~.$$
The cokernel of the symmetrization map in the Witt groups
$$\text{coker}(1+T:L_0(A) \to L^0(A))~=~\text{im}(L^0(A) \to 
Q_0(B(\infty),\beta(\infty)))$$ 
was computed for noetherian $A$ by Carlsson [5] in terms of 'Wu invariants'
prior to the general theory of Weiss [19].
\medskip

For $n\ge 0$ let
$(B\langle n+1 \rangle,\beta\langle n+1 \rangle)$ be the 
{\it $(v_{n+1}=0)$-universal chain bundle over $A$}, 
characterized up to equivalence by the properties 

{\parskip=2pt
\parindent=23pt
\item{(i)}
$\widehat{v}_r(\beta\langle n+1 \rangle):H_r(B\langle n+1 \rangle)
\to \widehat{H}^r(\Z_2;A)$ 
is an isomorphism for $r\ne n+1$,
\item{(ii)} $H_{n+1}(B\langle n+1 \rangle)=0$.\par}
\medskip

The {\it $(v_{n+1}=0)$-symmetric $L$-groups of $A$} are defined by
$$L\langle v_{n+1}\rangle^m(A)~=~
L\langle B\langle n+1 \rangle,\beta\langle n+1 \rangle\rangle^m(A)
~~ (m\ge 0)~.$$
Define the $A$-group 
$$Q\langle v_{n+1}\rangle~=~
Q\langle B\langle n+1 \rangle,\beta\langle n+1 \rangle\rangle^{2n}(S^nA)~,$$
to fit into the commutative braid of exact sequences 

$$\xymatrix@!C@C-30pt{
0 \ar[dr]\ar@/^2pc/[rr]^{} && 
Q_{2n}(S^nA) \ar[dr]^-{\di{r}} \ar@/^2pc/[rr]^-{\di{1+T}}&&
Q^{2n}(S^nA) \ar[dr]^-{\di{J}} \ar@/^2pc/[rr]^{}&&
0 \\
& \hbox{$\widehat{Q}^{2n+1}(S^nA)$} \ar[ur]^-{\di{H}} \ar[dr]^-{\di{1}} &&
Q\langle v_{n+1}\rangle \ar[ur]^-{\di{s}} \ar[dr] && 
\hbox{$\widehat{Q}^{2n}(S^nA)$} \ar[ur] \ar[dr] & \\
0 \ar[ur]\ar@/_2pc/[rr]_{} && 
Q^{2n+1}(S^nA) 
\ar[ur]\ar@/_2pc/[rr]^{\di{0}} &&
\hbox{$\widehat{Q}^{2n}(S^nA)$} \ar[ur]^-{\di{1}} \ar@/_2pc/[rr]_{} && 0 
}$$

\medskip

\noindent
In \S11 below we shall make use of the surjections 
$$L\langle v_{n+1}\rangle^{2n}(A) \to W_{Q\langle v_{n+1}\rangle}(A)~;~
(C,\phi,\theta,f) \mapsto (H^n(C),\phi_0,v_n(\phi,\theta,f))~~(n\ge 0)$$ 
defined for a field $A$. 
\medskip

\noindent \S11. {\bf An example}
\medskip
 
As an illustration of the exact sequence of \S8
$$\dots \to L_n(A) \to 
L\langle B,\beta\rangle ^n(A) \to Q_n(B,\beta)
~\raise4pt\hbox{$\partial\atop \to$}~L_{n-1}(A) \to \dots$$
we compute the Witt groups $L^0(A)$, $L_0(A)$, 
$L\langle v_1\rangle ^0(A)$ for $A$ a perfect field of characteristic 2, 
without appealing to the theorem of Arf [1] on the classification
of quadratic forms over such $A$ (cf. Example 2.14 of Ranicki [16]). 
\medskip

For any field $A$ $\cases L^{2n}(A)\cr L_{2n}(A)\endcases$ 
is the Witt group of nonsingular
$\cases \text{$(-1)^n$-symmetric}\cr
\text{$(-1)^n$-quadratic}\endcases$ forms over $A$, and
$\cases L^{2n+1}(A)=0\cr
L_{2n+1}(A)=0\endcases$ $(n\ge 0)$ - see Ranicki [15] for details. 
\medskip

Let then $A$ be a perfect field of characteristic 2, so that squaring
defines an automorphism 
$$A \to A~ ;~ a \mapsto a^2~.$$
Let $A$ have the identity involution 
$$\raise4pt\hbox{$\overline{~}$}~:~A \to A~ ;~ a \mapsto \overline{a}~=~a~.$$
As an additive group
$$\widehat{H}^r(\Z_2;A)~=~A~~(r \in \Z)$$ 
with $A$ acting by 
$$A \times \widehat{H}^r(\Z_2;A) \to \widehat{H}^r(\Z_2;A)~;~
(a,x) \mapsto a^2x~,$$
and there is defined an isomorphism of $A$-modules 
$$A \to \widehat{H}^r(\Z_2;A)~ ;~ a \mapsto a^2 ~.$$
The chain bundle over $A$ $(B(\infty),\beta(\infty))$ defined by 
$$\eqalign{
&d_{B(\infty)}~=~0~ :~B(\infty)_r~=~A \to B(\infty)_{r-1}~=~A~,\cr
&\beta(\infty)_s~=~\cases 1 \cr 0 \endcases~:~B(\infty)_r~=~A \to 
B(\infty)^{-r-s}~=~A~\text{if}~\cases s=-2r\cr s \ne -2r\endcases}$$
is universal. The twisted quadratic groups of $(B(\infty),\beta(\infty))$
are given up to isomorphism by
$$\eqalign{
&Q_{2n}(B(\infty),\beta(\infty))~=~Q^{\bullet}(A)~,\cr
&Q_{2n+1}(B(\infty),\beta(\infty))~=~Q_{\bullet}(A)~,}$$
with the abelian groups
$Q^{\bullet}(A)$, $Q_{\bullet}(A)$ defined by 
$$\eqalign{&Q^{\bullet}(A)~=~\{a\in A\,\vert\, a+a^2=0\}~=~\Z_2~,\cr
&Q_{\bullet}(A)~=~A/\{b+b^2\,\vert\, b\in A\}~,}$$
and isomorphisms defined by
$$\eqalign{
&Q_{2n}(B(\infty),\beta(\infty)) \to Q^{\bullet}(A)~;~
(\phi,\theta) \mapsto \phi_0(1)(1)~,\cr
&\hskip50pt\phi_0~:~B(\infty)^n~=~A \to B(\infty)_n~=~A~,\cr
&Q_{2n+1}(B(\infty),\beta(\infty)) \to  Q_{\bullet}(A)~;~
(\phi,\theta) \mapsto \theta_{-1}(1)(1)~,\cr
&\hskip50pt
\theta_{-1}~:~B(\infty)^{n+1}~=~A \to B(\infty)_{n+1}~=~A~.}$$

A symmetric form over $A$ $(M,\lambda)$ is even if and only if 
$$\lambda(x,x)~=~0 \in A~~ (x\in M)~ .$$ 
A nonsingular even symmetric form over $A$ $(M,\lambda)$ is
hyperbolic, since for any $x\in M$ there exists $y\in M$ such that
$\lambda(x)(y)=1\in A$, so that a hyperbolic summand may be split off
$(M,\lambda)$
$$(M,\lambda)~=~(Ax\oplus Ay,\pmatrix 0 & 1 \cr 1 &0\endpmatrix)
\oplus (M',\lambda')~,$$ 
with 
$$\text{rank}_AM'~=~(\text{rank}_AM)-2~.$$
Thus 
$$L\langle v_0\rangle ^0(A)~=~W\langle v_0\rangle (A)~=~0$$ 
and the symmetrization maps 
$$1+T~:~
L_{2n}(A)~=~L_0(A) \to L\langle v_0\rangle ^0(A) \to L^{2n}(A)$$
are zero. It is now immediate from the exact sequence 
$$\dots \to L_m(A)~\raise4pt\hbox{$1+T\atop \to$}~L^m(A) \to 
Q_m(B(\infty),\beta(\infty))~
\raise4pt\hbox{$\partial\atop\to$}~L_{m-1}(A) \to \dots$$
that 
$$L^{2n}(A)~=~Q^{\bullet}(A)~~,~~L_{2n}(A)~=~Q_{\bullet}(A)~.$$
In the symmetric case there is defined an isomorphism 
$$L^{2n}(A) \to Q^{\bullet}(A)~=~\Z_2~;~ (C,\phi) \mapsto 
\phi_0(v)(v)~=~\text{rank}_AH^n(C)~,$$ 
sending a $2n$-dimensional symmetric Poincar\'e complex $(C,\phi)$ over 
$A$ to the element $\phi_0(v)(v)\in Q^{\bullet}(A)$, with $v\in H^n(C)$ 
the unique cohomology class such that
$$\phi_0(x)(v)~=~\phi_0(x)(x) \in \widehat{H}^n(\Z_2;A)~~(x\in H^n(C))~.$$
The inverse isomorphism 
$$Q^{\bullet}(A)~=~\Z_2 \to L^{2n}(A)$$ 
sends $1\in Q^{\bullet}(A)$ to the $2n$-dimensional symmetric 
Poincar\'e complex $(C,\phi)$ defined by 
$$\phi_0~=~1~ :~ C^n~=~A \to C_n~=~A~ ,~ C_r~=~0~~ (r \ne 0)~.$$
In the quadratic case there is defined an isomorphism 
$$\partial~:~
Q_{\bullet}(A) \to  L_{2n}(A)~ ;~ a \mapsto (C,\psi)$$
with $(C,\psi)$ the $2n$-dimensional quadratic Poincar\'e complex over $A$ 
given by
$$\psi_0~=~\pmatrix a&1 \cr 0&1\endpmatrix~:~C^n~=~A\oplus A \to
C_n~=~A\oplus A~~,~~C_r~=~0~~ (r \ne n)~.$$
The inverse isomorphism $L_{2n}(A)\to Q_{\bullet}(A)$ sends a 
$2n$-dimensional quadratic Poinc\-ar\'e
complex $(C,\psi)$ over $A$ to the Arf invariant $c\in Q_{\bullet}(A)$ 
of the
nonsingular quadratic form $(H^n(C),(1+T)\psi_0,v^n(\psi))$ over $A$, as
defined by 
$$c~=~\sum\limits^m_{i=1}v^n(\psi)(x_{2i})v^n(\psi)(x_{2i+1}) \in 
Q_{\bullet}(A)$$
with
$\{x_i\,\vert\, 1 \le i \le m\}$ any basis for $H^n(C)$ such that 
$$(1+T)\psi_0(x_i,x_j)~=~\cases 1&\text{if $(i,j)=(2r,2r+1)$ or $(2r+1,2r)$}\cr
0&\text{otherwise~.}\endcases$$

The chain bundle over $A$
$(B\langle v_{n+1}\rangle ,\beta\langle v_{n+1}\rangle )$ 
$(n\ge 0)$ defined by
$$\eqalign{
&B\langle v_{n+1}\rangle _r~=~
\cases A&\text{if $r\ne n+1$}\cr 0&\text{if $r=n+1$}\endcases~,\cr
&d~=~0~:~B\langle v_{n+1}\rangle _r \to B\langle v_{n+1}\rangle _{r-1}~,\cr
&\beta\langle v_{n+1}\rangle _s~=~\cases
1 \cr 0\endcases~ :~ B\langle v_{n+1}\rangle _r \to 
B\langle v_{n+1}\rangle^{-r-s}~\text{if}~
\cases s=-2r,r\ne n+1\cr \text{otherwise}\endcases}$$
is $(v_{n+1}=0)$-universal. Define symmetric bilinear pairings 
$$\eqalign{&\rho~:~A \times A \to A~ ;~ (a,b) \mapsto ab~,\cr
&\sigma~ :~ Q^{\bullet}(A) \times Q^{\bullet}(A) \to A~ ;~ (1,1) \mapsto 1}$$
such that
$$\eqalign{
&Q\langle B\langle v_{n+1}\rangle ,\beta\langle v_{n+1}\rangle\rangle^{2n}
(S^nA)~=~A\times _{\rho}A~,\cr
&Q_{2n}(B\langle v_{n+1}\rangle ,\beta\langle v_{n+1}\rangle )~=~
A\times _{\sigma} Q^{\bullet}(A)~,\cr
&Q_{2n+1}(B\langle v_{n+1}\rangle,\beta\langle v_{n+1} \rangle )~=~0~.}$$
Let
$$\eqalign{Q\langle v_1 \rangle ~&=~Q\langle B\langle v_1\rangle ,
\beta\langle v_1\rangle\rangle^0(A)\cr
&=~A\times_{\rho}A~.}$$
Given a nonsingular
$Q\langle v_1\rangle$-quadratic form $(M,\lambda,\mu)$ over $A$
there exist $v\in M$, $\psi\in \Hom_A(M,M^*)$ such that 
$$\eqalign{&\lambda(x,y)~=~\lambda(x,v)\lambda(y,v) +
\psi(x)(y) + \psi(y)(x) \in A~~(x,y\in M)\cr
&\mu~ :~ M \to Q\langle v_1\rangle ~=~A\times_{\rho}A~;~
x \mapsto (\psi(x)(x),\lambda(x,v))~.}$$
The morphism 
$$\eqalign{&L\langle v_1\rangle ^0(A)~=~W_{Q\langle v_1\rangle}(A)
 \to Q_0(B\langle v_1\rangle ,\beta\langle v_1\rangle )~;\cr
&(M,\lambda:M\times M\to A,\mu:M\to Q\langle v_1\rangle ) \mapsto 
\mu(v)~=~(\psi(v)(v),\lambda(v,v))}$$
fits into a short exact sequence
$$0 \to L_0(A) \to L\langle v_1\rangle ^0(A) \to 
Q_0(B\langle v_1\rangle ,\beta\langle v_1\rangle ) \to 0~ .$$ 
The injection
$$L_{2n}(A) \to L\langle v_{n+1}\rangle ^{2n}(A) \to 
W_{Q\langle v_1\rangle }(A)~=~L\langle v_1\rangle ^0(A)$$ 
sends the cobordism
class of a $2n$-dimensional quadratic Poincar\'e complex over $A$ $(C,\psi)$
to the Witt class of the nonsingular $Q\langle v_1\rangle$-quadratic form
$(H^n(C),(1+T)\psi_0,\allowbreak iv^n(\psi))$, with $i$ the canonical injection 
$$i~:~H_0(\Z_2;A,(-1)^n)~=~A \to 
Q\langle v_1\rangle ~=~A\times _{\rho}A~;~ a \mapsto (a,0)~ .$$ 
\indent In the special case $A=\Z_2$
$$\eqalign{&L_{2n}(\Z_2)~=~Q_{\bullet}(\Z_2)~=~\Z_2~~,~~
L^{2n}(\Z_2)~=~Q^{\bullet}(\Z_2)~=~\Z_2~,\cr
&Q\langle v_{n+1} \rangle ~=~\Z_4~~,~~Q_0(B\langle v_1 \rangle,\beta
\langle v_1 \rangle)~=~\Z_4~,}$$
with
$$L\langle v_{n+1} \rangle^{2n}(\Z_2)~=~W_{\Z_4}(\Z_2)~=~\Z_8$$
the Witt group of nonsingular $\Z_4$-valued quadratic forms over $\Z_2$.
See Weiss [19,\t \S11] for the the algebraic Poincar\'e bordism
interpretation of the work of Browder [2] and Brown [4] on the Kervaire
invariant and its generalization, which applies also to the work of
Milgram [11].

\goodbreak

\enddocument